%% file: v12.15-20jun10.tex
\documentclass[11pt,reqno]{amsart}
\usepackage{amsmath, amssymb}
\usepackage[all]{xypic}
\newcommand{\eeq}{\end{equation}}

\newtheorem{theorem}{Theorem}[section]

\newtheorem{defi}[theorem]{Definition}
\newtheorem{lemma}[theorem]{Lemma}
\newtheorem{prop}[theorem]{Proposition}

\def\slfrac#1#2{\hbox{\kern.1em %
 \raise.5ex\hbox{\the\scriptfont0 #1}\kern-.11em %
 /\kern-.15em\lower.25ex\hbox{\the\scriptfont0 #2}}}

\font\phvr=phvr at 10pt
\newcommand\he[1]{\mbox{\phvr  #1}}

\newcommand{\eqn}[1]{(\ref{#1})}

\newcommand{\pf}{\noindent{\bf Proof.~}}
\newcommand{\beql}[1]{\begin{equation}\label{#1}}
\newcommand{\bsq}{{\vrule height .9ex width .8ex depth -.1ex }}

\newcommand{\La}{\Lambda}

\newcommand{\ep}{\epsilon}
\newcommand{\eps}{\epsilon}

\newcommand{\ZZ}{{\mathbb Z}}
\newcommand{\RR}{{\mathbb R}}

\newcommand{\PP}{{\mathbb P}}

\newcommand{\CC}{{\mathbb C}}

\newcommand{\hzt}{{\hat{L}}}
\newcommand{\ggG}{{\gamma}}

\newcommand{\bA}{{\bf A}}

\newcommand{\bh}{{\bf h}}

\newcommand{\sE}{{\mathcal E}}
\newcommand{\sF}{{\mathcal F}}

\newcommand{\sM}{{\mathcal M}}
\newcommand{\sN}{{\mathcal N}}
\newcommand{\sR}{{\mathcal R}}
\newcommand{\sS}{{\mathcal S}}

\newcommand{\hT}{{\he T}}

\title{The Lerch Zeta Function  I. Zeta Integrals}

\author{Jeffrey C. Lagarias}
\thanks{The research of the first author was supported by NSF grant DMS-0500555
and DMS-0801029 and that of the second author by NSF grant DMS-0457574
and DMS-0801096.}

\address{Department of Mathematics, University of Michigan,
Ann Arbor, MI 48109-1043,USA}
\email{lagarias@umich.edu}

\author{Wen-Ching Winnie Li}
\address{Department of Mathematics, Pennsylvania State University,
University Park, PA 16802-8401,USA}
\email{ wli@math.psu.edu}

\date{June 20, 2010}
\begin{document}
\begin{abstract} 
This is the first of  four papers that study  algebraic and
analytic structures associated to the Lerch zeta function.
This paper studies "zeta integrals" associated to the Lerch zeta function
using test functions, and obtains  functional equations for them. Special cases include 
a pair of symmetrized  four-term functional equations 
for combinations of Lerch zeta functions,
found by A. Weil, for real parameters $(a,c)$ with $0< a, c< 1$.
It extends these functions to real $a$, and $c$, 
and studies limiting cases of these functions where at least one of $a$ and $c$
take the values $0$ or $1$.   A main feature is that
 as a function of three variables $(s,a, c)$ with $a$, $c$ being real variables, 
the Lerch zeta function has discontinuities at integer values of
$a$ and $c$. For fixed $s$, 
the function $\zeta(s,a,c)$ is discontinuous on part of the boundary of the
closed unit square in the $(a,c)$-variables, and  the
location and nature of these discontinuities  depend on
the real part $\Re(s)$ of $s$.  Analysis of this behavior is used to determine 
  membership of these functions in $L^p([0,1]^2, da\, dc)$
 for $1 \le p < \infty$, as  a function of  $\Re(s)$.  The paper also defines
 generalized Lerch zeta functions associated to the oscillator representation,
 and gives analogous four-term functional equations for them.
\end{abstract} 

\maketitle
\tableofcontents

\setlength{\baselineskip}{1.2\baselineskip}

%
%
%
%
%

\section{Introduction}
The {\em Lerch zeta function}
\beql{101}
\zeta (s,a,c) := \sum_{n=0}^\infty e^{2 \pi ina} (n+c)^{-s}
\eeq
was introduced by Lipschitz \cite{Li57} in 1857, for real
$a$ and $c$ with $c > 0$, see also Lipschitz \cite{Li89}.
It is named after Lerch \cite{Le87}, who showed in 1887 that for
$\Im (a) > 0$ and $0 < c < 1$  it satisfies the three-term  functional equation
\beql{102}
\zeta (1-s, a,c) = (2 \pi )^{-s} \Gamma (s) \left\{
e^{\frac{\pi is}{2}} e^{-2 \pi iac} \zeta (s, -c,a)~
+~ e^{- \frac{\pi is}{2}} e^{2 \pi ic (1-a)} \zeta (s,c,1-a) \right\} ~,
\eeq
a result which is  called {\em Lerch's transformation formula}.
cf. Erdelyi \cite[p. 29]{Er53}.
Special cases of this function are the case $a=0$ which gives the
{\em Hurwitz zeta function}
\beql{103}
\zeta (s,c) := \sum_{n=0}^\infty \frac{1}{(n+c)^s} 
\eeq
studied by Hurwitz \cite{Hu82},
and the case $c=1$ which gives $e^{-2\pi i a} F(a, s)$, with
\beql{104}
F(a,s) := \sum_{n=1}^\infty
\frac{e^{2 \pi i na}}{n^s}
\eeq
being  the {\em periodic zeta function} studied in Apostol \cite[p. 257]{Ap76}. The Riemann zeta function
$\zeta (s)$ occurs as the intersection
of these two special cases, i.e. when $a = 0$ and $c=1$.

The Lerch zeta function has been extensively studied, see
the books of Lauren\v{c}ikas and Garunk\v{s}tis \cite{LG02},
Srivastava and Choi \cite[Chap. 2]{SC01} and Kanemitsu and Terada \cite[Chaps. 3-5]{KT07}
for many analytic results. There have been many proofs given of
Lerch's transformation formula and of equivalent functional equations, see for example
Apostol \cite{Ap51},
Oberhettinger \cite{Ob56}, Mikolas \cite{Mi71}, Berndt \cite{Be72}
and Weil \cite{We76}.

This is the first in a series of  four papers that studies
algebraic and geometric structures attached
to Lerch zeta functions. 
We begin with a brief summary of results in this series of
papers, and  give a more  detailed description of the main results of 
this paper in \S2.
%
%
%
%
%

\subsection{Overview of papers}

 In  Part I  our starting point is a pair of  symmetric functional
equations for the Lerch zeta function, viewing $a$ and $c$ as real
variables,  first given by A. Weil \cite{We76} in 1976.
These involve the two functions
$$
L^{\pm}(s, a,c) = \zeta(s, a, c) \pm e^{-2 \pi i a} \zeta(s, 1-a, 1-c).
$$
and relate the values of these  functions at parameter values $(s, a, c)$ to values 
of the same function at $(1-s, 1-c, a)$.
We derive functional equations for general zeta integrals of Lerch type
incorporating a test function, first for values $(a,c)$ on the unit square  
and then extended to all real values of $(a, c)$. We then deduce results
for $L^{\pm}(s,a,c)$ and the Lerch zeta function
on the extended domain by expressing them in 
terms of zeta integrals.
 We then  determine 
 that the Lerch zeta function 
has discontinuities
 at integer values of $a$ and $c$, and that the nature of these discontinuities
 depends on the 
real part of $s$.  
These discontinuities  appear to be  an important
aspect of the behavior of this function.

In Part II (\cite{LL2}) we analytically continue the
Lerch zeta function $\zeta(s, a, c)$ in three complex variables to a
 maximal domain of holomorphy. We show it  analytically continues to   a multivalued
function of $(s, a,c)$ which becomes single-valued on the 
maximal abelian covering manifold of the manifold 
$$
\sM :=\{ (s, a, c) \in  \CC \times (\CC \smallsetminus \ZZ) \times
(\CC \smallsetminus \ZZ)\}. 
$$
We determine the monodomy functions describing the multivaluedness.
Positive integer values of the $c\,$-variable turn out to be removable singularities, leading to 
a slightly larger analytic continuation. 
The remaining values where either  $a$ or $c$ are integers correspond to 
``singular strata" omitted from this analytic continuation. 
It is  possible to  define versions of the
Lerch zeta function restricted to various ``singular strata," which are {\em not part} of this analytic continuation.
In particular $a=0$ is a singular value, and both the 
Hurwitz zeta function (at $a=0,$ $c \not\in \ZZ$) and the Riemann
zeta function  (at $a=0, c=1$) live on ``singular strata".
It remains an interesting open problem to  further determine the
 relation of these  (degenerate) ``singular strata" functions
to the analytic continuation above, by limiting procedures; this provides
an approach to get more information on the behavior of  ``degenerate
functions" living on the singular strata.  From this perspective the discontinuity
properties studied in part I 
reflect  one aspect of the behavior of 
such limiting procedures. 

 In Part III (\cite{LL3}) we
make the variable change $z= e^{2 \pi i a}$ and analytically continue the
resulting {\em Lerch transcendent} 
$$
\Phi(s, z, c) = \sum_{n=0}^{\infty}  \frac{z^n}{(n+c)^{s}}
$$ 
in three complex variables. This function now lives on a covering manifold of
$$
\sN := \{ (s, z, c) \in \CC \times \left(\PP^{1}(\CC) \smallsetminus \{ 0, 1, \infty\} \right) \times (\CC \smallsetminus \ZZ)\}.
$$
It satisfies a  linear PDE  
over  the base $\sN$ having polynomial coefficients.
This function is related to the polylogarithm under the specialization $c=1$; here 
the specialization  $c=1$ lies on a regular stratum.
  We  investigate special values of $(s, a,c)$ where the monodromy
degenerates, and determine the degeneration.
 Such degeneration occurs for all values $s \in \ZZ$, which
 emphasizes  the  special significance of integer values of $s$,
 some of which are  ``critical values" in the sense of arithmetic geometry.

In Part IV (\cite{LL4}) we introduce a family of
two-variable ``Hecke operators" acting on  the $(a,c)$-variables, for
which the Lerch zeta function is a simultaneous eigenfunction. 
In suitable function spaces these satisfy
the Hecke relations  
$$\hT_m \hT_n= \hT_n \hT_ m = \hT_{mn}.$$
associated to $GL(1)$.
We determine properties of the  Hecke operator  action
in both the``real variables" context of part I and 
the ``complex variables" context of parts II and III. 
In the real-variables context we prove a uniqueness result characterizing,
for each fixed $s \in \CC$,  a two-dimensional
vector space of simultaneous eigenfunctions of these operators, together
with some side conditions, 
which generalizes a theorem of Milnor \cite{Mi83}
for the Hurwitz zeta function. 
In  the  complex-variables context there is an additional structure,
consisting of an induced  action of Hecke operators
on the infinite-dimensional vector space spanned by the monodromy functions.
This monodromy vector space also carries an action of the fundamental
group of $\sM$ and we determine information on  the commutation relations of the 
two-variable Hecke
operators with these operators. 

%
%
%
%
%

\subsection{Background and Motivation}

Our  original 
objective in this work was to seek  a geometric or 
``dynamical"  interpretation for the Lerch zeta function. 
 This was motivated
 by observations of
 the first author made in  \cite{La99}.
 He  studied two-variable dynamical zeta functions $\zeta(z, s)$,
 which included  both an arithmetical variable $s$
and a dynamical variable $z$, and observed that for certain dynamical systems 
(those being uniformly expansive in a suitable sense) the resulting dynamical zeta
function satisfied  an extra functional relation $\zeta(z,s) = \zeta(zq^{-s}, 0).$ 
This observation provided a framework to
explain the following coincidence:   the zeta function for a function field over a finite field
has an arithmetical definition using
the $s$-variable and also a dynamical definition using the $z$-variable, 
and these  match  under the change of variable  $z= q^{-s}$. 
An analogous result is unknown in the number field case.
The  Lerch zeta function is  a natural multiparameter generalization
of the Riemann zeta function, and it is natural to investigate 
whether its extra variables  might have  a geometric or ``dynamical"  interpretation.
A second objective, arising in the course of the work, was the
possibility of finding a representation-theoretic interpretation of
the Lerch zeta function. This possibility was based on  the feature that the Lerch
zeta function satisfies a linear partial differential
equation in the $(s, a, c)$-variables  which might play the role
of a Laplacian, and because it is a
simultaneous eigenfunction of a family of 
``Hecke operators," mentioned in the overview above.

These four  papers taken together
find various extra
structures attached to the Lerch zeta function,
whose form suggests the existence of a  more comprehensive theory (or theories)  in which they  fit,
addressing the objectives above.
In particular, there appear to be
two distinct underlying contexts in which the Lerch zeta
function appears, corresponding to the 
  ``real variables'' structure
and the ``complex variables'' structure mentioned in the summary above. The real variables
structure, which restricts the variables $a$ and $c$ to
be real,  imposes a twisted-periodicity condition described in this paper. 
The first author has since found
that this  structure does have a representation-theoretic interpretation
cf.  \cite{La09a}.
The complex variables structure, in which $a$ and $c$ are
complex variables,  and which leads to multivalued functions in
these variables,  seems  ``geometric'' in nature.  
We think it  will be an important
problem to determine interconnecting relations between these two
structures.

The results in these papers are formulated in a classical language. 
It seems likely that some of them may be reformulated and extended 
in an adelic framework. Other results  seem suitable for reformulation
in a framework of $D$-modules over  complex manifolds. 
We hope to return to these topics in subsequent work.

 It would be interesting to study
 if the structures studied in these papers have bearing on
 questions around 
  the Riemann  hypothesis.
It is well known that the Riemann hypothesis   
fails  to hold for the Lerch zeta function, but 
 the possible validity of the
weaker Lindel\"{o}f hypothesis is an open question. 
  Garunk\v{s}tis
and Steuding \cite{GS02} (see also Garunk\v{s}tis \cite{Ga05})
have put forth the proposal that the Lindel\"{o}f hypothesis may hold for the 
Lerch zeta function,  for  all real parameters $(a,c)$.
Our results in part II  are compatible with this proposal, in
that they imply that 
the Lindel\"{o}f hypothesis, if it holds for the Lerch
zeta function, would then also hold for 
its multivalued analytic continuation,  i.e.
it would  then hold for all
branches of the function lying over real values $(a,c)$.
Concerning   the Riemann hypothesis itself, 
the Lerch zeta function gives  a multi-parameter deformation of the
Riemann zeta  function,
in which it is located in a ``singular stratum" at a (non-isolated) singular point.
One may ask whether new information on the
Riemann hypothesis is obtainable through
taking various  degenerations approaching this ``singular stratum,"
using the extra structures available.

We conclude this discussion
by noting  that there  are 
further generalizations of the Lerch zeta function twisted by 
 Dirichlet characters $\chi~(\bmod~N)$. These take the form
 $$
 L^{\pm}(\chi, s ,a ,c ) := \sum_{n \in \ZZ} \chi(n) ({\rm sgn} (n+c))^k e^{2\pi i na} (n+c)^{-s},
 $$
 with $k=0,1$ and $\pm := (-1)^k$.
Our results can be extended to  apply to these
functions, since they can be expressed as linear combinations of scaled
versions of $L^{\pm}(s, a,c).$  We do not treat them in these papers, 
 in order to reduce the notational burden.\\

\paragraph{\bf Acknowledgments.} We thank K. Prasanna and P. Sarnak for helpful
comments. This  work began when the first
author was at AT \&T Labs and the second author visited there.
Both authors thank AT\&T Labs for support.

%
%
%
%
%

\section{Main results}
\setcounter{equation}{0}

 This paper proves variants of the functional equation of the Lerch zeta function,
in its four-term symmetrized form,  expressed using  the functions $L^{\pm}(s, a, c)$ introduced below,
and extends these results to zeta integrals for other test functions. These include test
functions for the oscillator representation.
Our focus is then on using these zeta integrals to  obtain information on  the
dependence of the Lerch zeta function on the $a$ and $c$
variables. In particular we determine continuity properties in
these variables as the variables approach integer values,
and membership of these functions in function spaces $L^p( [0,1]^2, da dc)$ for
various values of $s$.

The original functional equation given by Lerch \cite{Le87}
in 1887 was a three-term non-symmetric functional equation,
now called {\em Lerch's transformation formula. }
In 1976 A. Weil \cite[p. 57]{We76} gave 
 a pair of symmetrized four-term 
functional equations, which are the ones we consider here.
These two functional equations
encode an invariance under the action of the additive Fourier
transform on the real line.

To state  the symmetrized functional
equations  we introduce the two functions 
\beql{109b} 
L^{\pm}( s, a, c):=\zeta (s,a,c) \pm e^{-2 \pi ia} \zeta (s, 1-a, 1-c), 
\eeq
and we initially suppose the domain  is $0 < a, c < 1$. For
$\Re(s)> 1$ these functions  have the absolutely convergent Dirichlet
series representations
\begin{eqnarray}
L^{+}( s,a,c)  & = &
\sum_{n \in \ZZ \atop n \neq -c}
e^{2 \pi ina} |n+c|^{-s}, \\
L^{-}( s,a,c) & = &
 \sum_{n \in \ZZ \atop n \neq -c} sgn (n+c) e^{2 \pi ina} |n+c|^{-s} ~,
\end{eqnarray}
see Lemma~\ref{Jle22} below.

The Weil form of
the  functional equations involves the completion of these functions
obtained by adding appropriate  gamma factors.

\begin{theorem}\label{th11}
{\rm (Lerch Functional Equations).} Let $a, c$ be real with $\{ (a, c) \in (0,1) \times (0,1)\}.$ 

(1) The completed  function
\beql{110a}
\hat{L}^{+}(s,a,c)   :=\pi^{- \frac {s}{2}} \Gamma (\frac {s}{2}) L^{+}(s, a, c)
\eeq
analytically continues to
 an entire function  of $s$, and
satisfies the functional equation
\beql{110b}
 \hat{L}^{+}( s,a,c)= e^{- 2 \pi ia c} \hat{L}^{+}( 1-s, 1-c, a).
\eeq

(2) The completed function
\beql{111a}
\hat{L}^{-}( s,a,c) := \pi^{- \frac {s+1}{2}} \Gamma (\frac {s+1}{2}) L^{-}( s,a,c)
\eeq
analytically continues to
 an entire function of $s$, and satisfies the functional equation
\beql{111b}
 \hat{L}^{-}(s,a,c)= i e^{- 2 \pi ia c} \hat{L}^{-}( 1-s, 1-c, a).
 \eeq
\end{theorem}

The fact that 
the ``completed" functions
$\hat{L}^{\pm}(s, a, c)$ are entire functions of $s$ for fixed $0<a, c<1$
implies that for these values the (non-completed)
function $L^{+}(s, a, c)$ necessarily has ``trivial zeros" at $s= 0, -2, -4, ...$, while
the (non-completed) function $L^{-}(s, a,c)$ has ``trivial zeros" at $s=-1, -3, -5, ...$,
for all values of $a$ and $c$.

In \S3 we derive these functional equations using a method analogous to 
 Tate's thesis \cite{Ta67}, assuming $a, c$ are real with $0< a, c <1.$  We consider a
{\em   zeta integral of Lerch type}  $F_k(f; s, a, c)$ 
attached to  a general test function $f(x)$ in
the Schwartz space $\sS(\RR)$, with $k=0, 1$,
and analytically continue  these integrals to  entire functions of $s$,
for fixed $(a, c)$ as above.
The perspective of Tate's thesis has two important features
 which are: (i) to realize a ``zeta function''
as a greatest common divisor of a set of ``zeta integrals'' indexed
by test functions; and (ii) to carry this out both in local settings and in a global adelic
setting. (Compare Ramakrishnan and Valenza \cite[p. 242]{RV99}.)
Here we  carry out
 for Lerch zeta functions 
the ``zeta integral'' step using a set of test
functions at the real place which transform nicely under
the Fourier transform.  
Our treatment is global rather than local in
that it uses Poisson summation, and it emphasizes
 operator aspects
 in the use of ``averaging" operators. 
An adelic treatment  of the Lerch zeta function
 involves further issues, which we do not take up here.

In \S4 we extend the variables in   Lerch zeta integrals  $F_k(f; s, a, c)$,
to all values $(a, c) \in \RR \times \RR$.
We derive a general zeta integral functional equation
 (Theorem~\ref{JLth21}) for $F_k(f; s, a, c)$,
specified by a test function $f(x)$ in
the Schwartz space $\sS(\RR)$, with $k=0, 1$.
 For integer values of $a$ or $c$
the functions are no longer entire functions, but are meromorphic
functions, having  possible simple poles  at $s=0$ and $s=1$. 
We observe that these functions exhibit   extra ``twisted periodicity" relations in the $(a, c)$-variables. 

In \S5, on  specializing  to the Gaussian test functions $\phi_0(x)= e^{- \pi x^2}$
and $\phi_1(x) = xe^{\pi x^2}$, we  obtain
extensions to all $(a,c) \in \RR \times\RR$ 
of the functions $L^{\pm}(s, a,c)$ and from these define  an extension
$\zeta_{\ast}(s, a, c)$ of the Lerch zeta function to $(a,c) \in \RR \times\RR$ 
which preserves the symmetrized functional equations. 
The extended function  $\zeta_{\ast}(s, a,c)$  is an entire function of
$s$ whenever both  $a$ and $c$ are not integers,  and otherwise 
it is meromorphic in $s$, with
possible simple poles located only at $s=0$ or $s=1$. 
We show that for $\Re(s) > 1$ it
is explicitly given for fixed $(a, c)$ by 
\beql{sp111b} 
\zeta_{\ast}(s,a, c) := \sum_{n + c > 0} e^{2 \pi i n a} (n + c)^{ -s}, 
\eeq 
 The extended functional
equations are as follows.

\begin{theorem}\label{th12}
{\rm (Extended Lerch Functional Equations)}
 Let $(a,c) \in \RR \times \RR$. The extended Lerch zeta function $\zeta_{*}(s, a, c)$ is
a meromorphic function of $s \in \CC$,
which  has a simple pole at
$s = 1$ if $a \in \ZZ$, and is holomorphic in $\CC$ otherwise.

(i) In the $a$-variable it is periodic, with
\beql{sp111c}
\zeta_{*}(s, a + 1 , c) = ~~~\zeta_{*}(s, a, c).
\eeq

(ii) In the $c$-variable it
satisfies  the twisted periodicity equation
\beql{sp111d}
\zeta_{*}(s, a, c+1) = e^{ - 2 \pi i a} \zeta_{*}(s, a, c).
\eeq

(iii) For $k=0$ and $1$ and $\pm = (-1)^k$ the completed extended  functions  
\beql{sp111e}
\hat{L}_{\ast}^{\pm}( s,a,c) := \pi^{- \frac{s+k}{2}} \Gamma \left(
\frac{s+k}{2} \right)  \left( \zeta_{\ast} (s,a,c) + (-1)^k e^{- 2 \pi ia} \zeta_{\ast} (s,1 -a , 1 - c) \right) 
 \eeq
 are meromorphic for $s\in \CC$, and are analytic except for possible  simple poles at $s = 0$
  or $1$. The poles occur only for  $k=0$, and then  if and only if $a \in \ZZ$
or  $c \in \ZZ$.
The completed extended functions 
satisfy the functional equations, for $k=0,1$, 
\beql{sp111f}
\hat{L}_{\ast}^{\pm}(s,a,c) = i^k e^{-2 \pi iac}\hat{L}_{\ast}^{\pm}(1-s, 1-c,a).
\eeq
\end{theorem}

The extended functions $L_{\ast}^{\pm}(s,a,c)$ inherit the same symmetries as
$\zeta_{\ast}(s, a, c)$; namely,  for all $(a,c) \in \RR \times \RR$,
the twisted periodicity conditions hold:
\begin{eqnarray*}
 L_{\ast}^{\pm}( s, a + 1 , c) &= &~~~~ L_{\ast}^{\pm}( s, a, c), \\
 L_{\ast}^{\pm}(  s, a, c+1) & = &
e^{ - 2 \pi i a}L_{\ast}^{\pm}( s, a, c).
\end{eqnarray*}

Theorem \ref{th12}  defines the Lerch zeta function $\zeta_{\ast}(s, a, c)$
for each $(a,c)$ on the boundary of the unit square, for all $s\in \CC$, 
except for $s=0, 1$ where poles occur.  Using \eqn{sp111b} and the
functional equations above we find that on the boundary it is
given in terms of the Hurwitz zeta function and periodic zeta
function by 
 \beql{112a}
  \zeta_{\ast}(s, 0, c)= \zeta(s, c) ~~~~~~~~~~~~~for~~0 < c \leq 1,
\eeq
\beql{112b}
\zeta_{\ast}(s, 1, c) = \zeta(s, c) ~~~~~~~~~~~~~for~~0 < c \leq 1,
\eeq
\beql{112c} \zeta_{\ast}(s, a, 0) = F(a, s) ~~~~~~~~~~~~for~~0
\leq a \leq  1,
\eeq
and
\beql{112d}
 \zeta_{\ast}(s, a, 1) = e^{- 2\pi i a} F(a, s) ~~~~for~~0 \leq a \leq 1.
\eeq 
  At the four corners
of the square we have 
\beql{112e}
 \zeta_{\ast}(s, 0, 0) = \zeta_{\ast}(s, 0,1) = \zeta_{\ast}(s, 1, 0) = \zeta_{\ast}(s, 1, 1)= \zeta(s),
\eeq
 the Riemann zeta function. Theorem~\ref{th12} therefore includes
meromorphic continuations of the Hurwitz zeta function and the
periodic zeta function to all $s \in \CC$, and also gives functional
equations for these functions. 
At the end of \S5, we deduce  an extended  Lerch
transformation formula valid for all $(a,c) \in \RR \times \RR$
by taking appropriate combinations of
 $\hat{L}_{\ast}^{+}( s, a, c)$ and $\hat{L}_{\ast}^{-}( s, a, c)$
and using suitable
Gamma function identities (Theorem~\ref{th51}).

In \S6 we analyze how the extended Lerch zeta function $\zeta_{\ast}(s, a,c)$
behaves  as $a$ and $c$ approach integer values.
We observe that Theorem~\ref{th12} comes at a price: the function
 $\zeta_{\ast} (s, a, c)$ is {\em discontinuous at integer values of  the $c$ and $a$
variables, for certain ranges of $s$.} Indeed the formula 
\eqn{sp111b} reveals that  for $\Re(s) > 1$ the function $\zeta_{\ast}(s,a, c)$ is  
 discontinuous at integer values of the $c$-variable. 
 It turns out  to be discontinuous at integer values of the $a$-variable
 for another range of $s$-values. 
These ranges of $s$ differ for the $a$-variable and the
$c$-variable. 
Let $\Box = \{ (a, c): ~ 0 \le a \le 1,~0 \le c \le1\}.$
We establish a  limiting formula for
$\zeta_{\ast} (s,a,c)$ with
$(a,c) \in \Box^\circ$ as it approaches the boundary of $\Box$,
which exhibits the discontinuities (Theorem \ref{th31}).
Using it, we obtain  precise conditions when
$\zeta_{\ast}  (s,a,c)$ continuously extends to portions of the boundary
of the unit square. Note  that the equality $\zeta(s, a, c)= \zeta_{\ast}(s, a,c)$ holds
for $(a,c) \in \Box^\circ$, so we may state these results for $\zeta(s,a,c)$, as follows.

\begin{theorem}\label{th14}
{\rm (Continuous Extension to Boundary)}
For fixed $s \in \CC$, the function
$\{ \zeta (s,a,c) : (a,c) \in \Box^\circ \}$ continuously
extends to portions of the boundary $\partial \Box$
of the unit square to define a function $\zeta (s, a, c)$ there,
in the following cases:
\begin{itemize}
\item[(i)]
$a = 0$ and $0 < c < 1$ when $\Re (s) > 1;$
\item[(ii)]
$a = 1$ and $ 0 < c < 1$ when $\Re (s) > 1;$
\item[(iii)]
$c = 0$ and $ 0 < a < 1$ when $\Re (s) < 0;$
\item[(iv)]
 $c = 1$ and $0 < a < 1$  when $ s \in \CC;$
\item [(v)]
the two corners $(a,c) = (0,1)$ and $(1, 1)$, when $\Re (s) > 1.$
\end{itemize}
In these cases this continuous extension satisfies
\beql{sp113}
\zeta (s, a, c) = \zeta_{*} (s, a, c).
\eeq
A continuous extension to the boundary does not exist for any other
values of $(s,a,c)$.
\end{theorem}

\begin{figure}[htb]
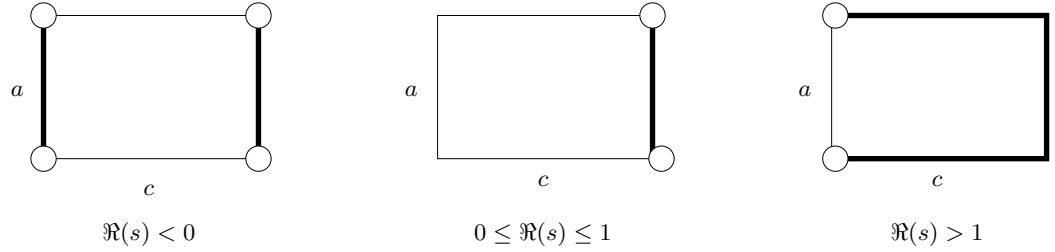

\begin{center}
\input fg1.pstex_t
\end{center}
\caption{Continuous extension of $\zeta (s,a,c)$ (heavy lines).}
\label{fg101}
\end{figure}

The ranges where continuous extensions are possible are
indicated schematically by the heavy lines in Figure~\ref{fg101}.
Interesting
features to note about the discontinuities are:

(1) The locations
of the discontinuities
{\em depend only on the value of $\Re(s)$.}

(2) The locations of the discontinuities are
{\em not invariant under the transformation
$(s, a, c)  \rightarrow (1 - s, 1 - c, a)$} that appears in the
functional equation.

(3) There is no continuous extension to the corner points
$(a,c) = (0,0)$ and $(1,0)$. 

The occurrence and behavior of the discontinuities is an important
feature of the Lerch zeta function. 
In part II we show  that the integer values of
$a$ and $c$ correspond to complex singularities to these functions, 
on  viewing
$a$ and $c$ as complex variables. From this perspective these discontinuities
form part of the  structure of these singularities.  It is 
important to understand the nature of these singularities because  the 
Hurwitz zeta function, Periodic zeta function and Riemann zeta
function all arise in the cases where one or both $a$ and $c$ are integers.

In  \S7 we use the determination of the limiting behavior of
$\zeta_{\ast}(s, a,c)$ to define ``renormalized'' versions $L^{R, \pm}(s, a,c)$
of the Lerch function $L^{\pm}(s,a,c)$ that continuously extend to
the boundary of $\Box$, for $s \in \CC  \smallsetminus \ZZ$.
The precise result is given as
Theorem~\ref{th31A}. The ``renormalized" functions 
remove four terms which account for the
divergent behavior at the singularity.

More generally, for fixed $s \in \CC$, let $\sE_s$ (for ``Eigenspace")
denote the complex vector space  spanned by the four functions
\beql{N241}
\sE_{s} := \CC[ L^{\pm}(s, a, c), e^{-2\pi i ac}L^{\pm}( 1-s, 1-c, a)].
\eeq
The  functional equation gives linear  dependencies
among these four functions, which show that, for each $s \in \CC$,
the vector space $\sE_s$ is  two-dimensional. In fact we have
\begin{eqnarray} 
\sE_{s} &=& \CC[ L^{\pm}(s, a, c)]   ~~~~~~~~~~\mbox{for}~~ s \in \CC \smallsetminus \ZZ_{\le 0}
\label{242} \\
\sE_{s} & = & \CC[ e^{-2 \pi i  a c}L^{\pm}(1-s, 1-c, a)]  ~~\mbox{for}~~ 
s \in \CC \smallsetminus \ZZ_{\ge 1}.
\label{243}
\end{eqnarray}
Note that 
$
\zeta(s, a, c) = \frac{1}{2} \left( L^{+}(s, a,c)+ L^{-}(s, a,c) \right)$ belongs to $\sE_s.
$
We can now ``renormalize" all functions in $\sE_s$ using the renormalized
versions of the basis functions. \\

We use the  renormalized Lerch zeta functions to  determine membership of all
elements of $\sE_s$  in   the function space $L^p( \Box, da~ dc)$ for $1 \le p \le \infty$,
as follows.

\begin{theorem}\label{th15} {\rm ($L^p$-Membership)}

(1) Suppose $1 \le p < 2$. Then all functions in $\sE_s$
belong to  $L^p(\Box, da ~dc)$ if and only if $s \in \CC$ satisfies  
$$1 - \frac{1}{p}< \Re(s) < \frac{1}{p}.$$
Otherwise, 
no nonzero function in $\sE_s$ belongs to 
$L^p(\Box, da~ dc)$.

(2)  Suppose $2 \le p \le \infty$. Then
no nonzero function in $\sE_s$ belongs to 
 $L^p(\Box, da ~dc)$, for all $s \in \CC$.
\end{theorem}

The most interesting cases are $p=1$, where 
$L^{\pm}(s,a ,c)$ both belong to  the Banach space
$L^1(\Box, da~ dc)$ for $s$ inside the critical strip  $0 < \Re(s) <1$,
and $p=2$, where no nonzero functions in  $\sE_s$ belong to the
Hilbert space $L^2(\Box, da ~dc)$, for all $s \in \CC$. 
Theorem~\ref{th15} will be relevant to the Hilbert spaces studied in Part IV.

In \S8 we  show  compatibility of the Lerch
zeta function with the oscillator representation.
The Schwartz space $\sS(\RR)$ is the set of smooth
vectors of a representation of (a double cover of) the group $SL(2, \RR)$,
the {\em oscillator representation},
and also of the Heisenberg group acting on $L^2(\RR, dx)$, a
viewpoint emphasized in Howe~\cite{Ho80}. It
is closed under the action of the Fourier transform, and also closed under the
action of the Weyl algebra $\bA_1 = \CC [ x, \frac{\partial}{ \partial{x}}]$.
The Mellin transform acts as an intertwining operator,
converting the Weyl algebra action into  difference operators in the $s$-variable,
and this  action is also compatible with the Fourier transform.
Bump and Ng \cite{BN86} and Bump et al \cite{BCKV99} observed
for the Riemann zeta function
that there is an an infinite
family of functional equations related to the oscillator
Hamiltonian. We show that this phenomenon generalizes
to the Lerch zeta function, with a similar proof.  We 
 obtain an infinite family of
functional equations for functions $L_n(s, a,c)$ associated to the eigenfunctions of the
oscillator Hamiltonian, indexed by the nonnegative
integers $n \ge 0$. 
 The functional equations in
Theorem~\ref{th11} correspond to the cases $n=0$ and $n=1$, respectively,
whose test functions $\phi_0(x)$ and $\phi_1(x)$ are vacuum vectors for the two 
irreducible subrepresentations of
the oscillator representation. The result of Bump and Ng \cite{BN86} for  the Riemann zeta function
is  then in principle recoverable by specialization of variables. \\

\paragraph{\bf Notation.}
The hat notation
$\hat{L}$ always denotes a ``completed''  function multiplied
by an appropriate archimedean Euler factor (gamma factor).
For functions $f$ on the real line
we  use $\sF f$ to denote an additive Fourier transform
and $\sM f$ to denote a Mellin transform (multiplicative Fourier
transform), as defined in \S2.
The variable $\epsilon \in \{ -1, 1\}$ and we often write
$\epsilon= (-1)^k$ for $k = 0, 1$. The function $\delta_{\ZZ}(x)=1$
if $x \in \ZZ$ and is $0$ otherwise. Here  $\Re(s)$ and $\Im(s)$ denote
the real and imaginary parts of a complex variable $s$.

%
%
%
%
%
\section{Analytic continuation for zeta integrals}
\setcounter{equation}{0}

The Schwartz space $\sS (\RR )$ consists of all
smooth functions $f(x)$ such that it and all its derivatives decrease 
more rapidly than any power of $|x|$ as $x \to \pm \infty$.
We decompose it into  even
and odd functions
\beql{201}
\sS (\RR) = \sS^+ \oplus \sS^-,
\eeq
and write
\beql{202}
f(x) = f_+(x) + f_-(x),
\eeq
with
$$
f_+(x) := \frac {1}{2}(f(x) + f(-x)) \quad\mbox{and}\quad
f_-(x) := \frac {1}{2}(f(x) - f(-x)).
$$
The {\em Fourier transform} $\sF f(y)$ is given by \beql{203} \sF
f(y) := \int_{-\infty}^\infty e^{ - 2 \pi i xy} f(x) dx. \eeq It
takes $\sS^+$ to  $\sS^+$, and  $\sS^-$ to $\sS^-$, so that
\beql{203a}
 \sF f_{\epsilon} = (\sF f)_\epsilon
\qquad\mbox{for}\qquad \epsilon= (-1)^k,
 \eeq
  and also  satisfies
$\sF \sF f(x) = f(-x)$. The functions $\phi_0(x) = e^{-\pi x^2}$ and
$\phi_1(x) = x e^{-\pi x^2}$ have Fourier transforms \beql{eq204}
\sF \phi_k(y) = (-i)^k \phi_k (y)~~~~for~~~k=~0,~1. \eeq

We use the {\em two-sided Mellin transforms} $\sM_k$ with $k = 0,
1~(\bmod~ 2)$, defined by 
\beql{205}
 \sM_k (f)(s) :=
\int_{-\infty}^\infty f(x)(sgn~ x)^k |x|^s \frac{dx}{|x|} ~~~ for~~~k=~0,~1. 
\eeq 
The {\em (one-sided)
Mellin transform} $\sM$ on $\RR^+$ is given by 
 \beql{206}
  \sM (f)(s) := \frac{1}{2} (\sM_0 (f)(s) + \sM_1(f)(s)) = \int_{0}^\infty f(x) x^s \frac{dx}{x}.
 \eeq
 We need only consider $\sM_0$ acting on even
Schwartz functions and $\sM_1$ acting on odd Schwartz functions,
namely 
\beql{206a} \sM_k(f)(s) = \sM_k (f_\epsilon)(s)
\quad\mbox{where}\quad \epsilon = (-1)^k,
\eeq 
since we have
 \beql{207} 
 \sM_k( f_\epsilon)(s) \equiv 0
\quad\mbox{when}\quad \epsilon ~=~ (-1)^{k+1}~.
 \eeq

The functional equations are associated  to the action of certain
multiplicative averaging operators on test function $f \in
\sS(\RR)$. \\

\begin{defi}~\label{nde20}
{\em
(1) For fixed $0 < a , ~ c < 1$ , and $f (x) \in \sS(\RR)$  the 
{(\em multiplciative) averaging
operator } $A^{a,c}$ is given by
\beql{209}
 A^{a,c}[f](x) := \sum_{n \in \ZZ}
f((n+c)x)  e^{2 \pi ina} ~. 
\eeq 
The rapid decrease of $f(x)$ at
$\pm \infty$ ensures that $A^{a,c}[f](x)$ is well-defined and
continuous on $\RR \smallsetminus \{0\}$, however it may have a
singularity at $x=0$. Its Mellin transform $\sM
(A^{a,c}[f])(s)$ is well-defined for $\Re (s)> 1$.

(2) The
{\em symmetrized (multiplicative)  averaging operators}
$A_k^{a,c}$ on $f \in \sS(\RR)$are given by
\beql{210b}
A_k^{a,c}[f](x):= 
A^{a,c}[f](x)  + (-1)^k e^{-2 \pi ia} A^{1-a, 1-c} [f](x).
\eeq
That is, 
\begin{eqnarray}\label{210}
A_k^{a,c}[f](x)
& = & \sum_{n \in \ZZ} f((n+c) x)  e^{2 \pi ina}  \nonumber \\
&& \quad + (-1)^k
e^{-2 \pi ia} \sum_{m \in \ZZ}   f((m+1-c)x) e^{- 2 \pi ima} ~.
\end{eqnarray}
}
\end{defi} 

Note that 
\beql{211}
A^{a,c}[f](x) := \frac{1}{2}(A_0^{a,c} + A_1^{a,c} )[f](x) = 
 \sum_{n \in \ZZ} f((n+c)x) e^{ 2\pi i n a}.
\eeq

Our main object of study will be the Mellin transforms of functions
acted on by these averaging operators.
\begin{defi}\label{de20}
{\rm
(1) For the test function $f \in \sS(\RR)$ and $k=0,1$, 
the  {\em zeta integral} $F_k (f; s,a,c)$  is
defined for  $0 < a<1$, $0<c<1$
to be the  one-sided Mellin transform of $A_{k}^{a,c}[f]$, phase-shifted
by a factor $e^{\pi i a c}$, i.e. 
namely
\beql{212}
F_k (f; s,a,c)  :=  e^{\pi iac} \sM (A_k^{a,c}[f])(s) =e^{\pi iac} \int_0^\infty A_k^{a,c}[f](x) x^{s-1} dx.
\eeq

(2) The {\em full zeta integral} $F(f; s , a, c)$ for the test function $f$ is 
\beql{eq213}
F(f; s,a,c) :=  2e^{\pi iac} \sM (A^{a,c}[f])(s) = F_0(f; s, a, c) +
F_1(f; s, a, c).
\eeq
}
\end{defi}

These Mellin transforms are related to the Lerch zeta function
as follows.
\begin{lemma}\label{nle21}
Let $0 < a < 1$ and $0 < c < 1$.
Then for $f \in \sS(\RR)$ and $\Re (s) > 1$ and $k=0,1$,
\beql{213a}
 \sM (A_k^{a,c}[f])(s)  =  \sM_k(f)(s)  L^{\pm}( s, a, c)
\eeq
in which  $\pm = (-1)^k$ and
\beql{213b}
 L^{\pm}( s, a, c)    : = \zeta (s,a,c) \pm  e^{-2 \pi ia} \zeta (s,1-a,1-c) ~.
\eeq
Thus 
\beql{213} 
F_k (f; s, a, c) =  e^{\pi i  a c} \sM_k(f)(s) L^{\pm}( s, a, c).
\eeq
\end{lemma}

\paragraph{\bf Remark.}
For $(a,c) \in \Box^\circ$ and $\Re (s) > 1$ we have
 \beql{214}
L^{\pm}(s,a,c) = \sum_{n \in \ZZ} e^{2 \pi ina} (sgn\,(n+c))^k  | n+c|^{-s}    ~.
\eeq
 Lemma \ref{nle21} may be viewed as
showing that $e^{\pi iac}L^{\pm}(s,a,c)$ is the ``symbol'' of
the averaging operator $B_k^{a,c}$ under the Mellin
transform.\\

\pf
We first set
$$
\Lambda_k (f; s,a,c) := \sM_k (f)(s)  \left(e^{\pi iac} \zeta (s,a,c) \right).
$$
The formula \eqn{213} is equivalent to showing that
\beql{215}
F_k (f;s,a,c) = \Lambda_k (f; s,a,c) +
(-1)^k e^{-\pi i(a+1-c)} \Lambda_k (f; s,1-a, 1-c) ~.
\eeq
To verify \eqn{215} we note that for $\Re (s) > 1$,
$$
e^{\pi iac} \zeta (s,a,c)= \sum_{n=0}^\infty e^{2 \pi ia (n+ c/2)} (n+c)^{-s} ~,
$$
hence
\begin{eqnarray}\label{JL29}
\La_k (f; s,a,c) & = &
\left( \sum_{n=0}^\infty e^{2 \pi ia (n+\frac {c}{2})} (n+c)^{-s} \right)
\int_{-\infty}^\infty
f(x) (sgn (x))^{k} |x|^{s - 1}~ dx  ~.
\nonumber \\
& = & \int_0^\infty \left( \sum_{n=0}^\infty e^{2 \pi ia (n+ \frac{c}{2})}
f((n+c)x)) \right) x^{s-1} dx \nonumber \\
&& + (-1)^k \int_0^\infty \left( \sum_{n=0}^\infty e^{2 \pi i (-a) (-n - \frac{c}{2} )}
f(-(n+c)x) \right) x^{s-1} dx.
\end{eqnarray}
Combined with a similar formula for $\Lambda_k (f; s,1-a,1-c)$, a
detailed calculation yields
\begin{eqnarray}\label{217}
\Lambda_k (f; s,a,c) &+& (-1)^k e^{- \pi i (a+1-c)} \Lambda_k (f;s,1-a,1-c)
\nonumber \\
&=& e^{2 \pi iac} \left\{
\int_0^\infty \sum_{n \in \ZZ} e^{2 \pi ina} f ((n+c)x) x^{s-1} dx \right.
\nonumber \\
&& 
+ (-1)^k
e^{-2 \pi ia} \left.
\int_0^\infty \sum_{m \in \ZZ}
e^{-2 \pi ima} f((m+1-c)x) x^{s-1} dx \right\} 
\end{eqnarray}
as required.~~~$\bsq$ \\

Lemma \ref{nle21} shows that the dependence of $\sM (A_k^{a,c}[f])(s)$ on
the test function $f$ is confined to the
Mellin transform $\sM_k (f)(s)$, hence we have
\beql{JL210a}
 F_k (f; s,a,c) = F_k (f_\epsilon ; s,a,c)
\quad\mbox{where} \quad \epsilon = (-1)^k.
\eeq
Substituting \eqn{JL29} into this definition yields the
integral representation
\begin{eqnarray}\label{JL211}
F_k (f; s,a,c) ~~~~~ = & \int_0^\infty \left( \sum_{n \in \ZZ}
e^{2 \pi ia (n+ \frac{c}{2} )}
f((n+c)x) \right) x^{s-1} dx   \\
 +  (-1)^ke^{- \pi i(a+1-c)}&
\int_0^\infty \left( \sum_{m \in \ZZ} e^{2 \pi i (1-a) (m+
\frac{1-c}{2} )} f((m+1-c) x) x^{s-1} dx \right). \nonumber
\end{eqnarray}
For even or odd Schwartz functions we have
 \beql{JL212} 
 F_k(f_{\ep};  s,a,c) = \left\{
\begin{array}{lll}
2 \int_0^\infty \left( \sum_{n \in \ZZ}
e^{2 \pi ia (n+c/2)} f((n+c)x) \right)
x^{s-1} dx & \mbox{if} &
\ep = (-1)^k, \\ [+.2in]
0 & \mbox{if} & \ep = (-1)^{k+1}.
\end{array}
\right.
\eeq
This follows by pairing the terms $n$ and $m=-n-1$ in the two integrals
in \eqn{JL211}:
they either cancel or match according to the value of $k$.

We use the Poisson summation formula to obtain an analytic
continuation and functional equation for $F_k (f; s,a,c)$, and, in \S4, to 
extend the definition to $(a,c) \in \RR \times \RR$. This involves
splitting the integral representation \eqn{JL211} into two pieces,
 $\int_0^1$ and $\int_1^\infty$. We give
the $\int_1^\infty$ pieces the name $\Phi_k (f; s,a,c)$. If we
define $\Phi (f; s,a,c)$ by
 \beql{JL213b} \Phi (f; s,a,c) := \int_1^\infty \left( \sum_{n \in \ZZ} e^{2 \pi ia (n+ c/2)}f((n+c)x)
\right)x^{s-1} dx.
 \eeq 
 then we have, for $k=0, 1$,
\beql{JL213a}
 \Phi_k (f; s,a,c) :=  \Phi (f; s,a,c) +
(-1)^k e^{- \pi i (a+1-c)} \Phi (f; s,1 - a, 1 -c).
\eeq
The function $\Phi (f; s,a,c)$ defined by \eqn{JL213b} is an entire
function of $s$ on the domain $0<a<1, 0<c<1$,
and the same property is inherited by $\Phi_k (f; s,a,c)$
for $k = 0,1.$

We now show:
\begin{lemma}\label{le21}
For fixed $0 < a < 1$ and $0 < c < 1$ and $f(x) \in \sS(\RR)$, the
function $F_k (f; s,a,c)$ is given by 
\beql{JL214}
 F_k (f; s,a,c)=
\Phi_k (f; s,a,c) + (-1)^k e^{- \pi ia} \Phi_k ( \sF f;  1-s, 1-c,a)
 \eeq
 for $k=0$ or $1$ and $\Re (s) > 1.$ The right side
analytically continues $F_k (f; s,a,c)$ to an entire function of $s
\in \CC$.
\end{lemma}

\pf
The integral representation \eqn{JL211} gives
\beql{JL215}
F_k (f; s,a,c) = \Phi_k (f; s,a,c) + \Psi_k (f; s,a,c),
\eeq
in which $\Phi_k(f; s, a,c)$ is given by \eqn{JL213a} and 
\begin{eqnarray}\label{JL216}
&&\Psi_k (f; s,a,c) := \int_0^1 \left( \sum_{n \in \ZZ} e^{2 \pi
ia (n+ \frac{c}{2} )}
f((n+c)x) \right) x^{s-1} dx   
\nonumber\\
&&+ (-1)^k e^{- \pi i (a+1-c)}
\int_0^1 \left( \sum_{n \in \ZZ} e^{2 \pi i (1-a) (n+ \frac{1-c}{2} )}
f((n+1-c) x) \right) x^{s-1} dx \,. 
\end{eqnarray}
Write $f= f_+ + f_-$ as a sum of even and odd functions, and we obtain 
 \beql{JL217}
  \Psi_k ( f_\ep ; s,a,c) =  \left\{
\begin{array}{lll}
2 \int_0^1 \left( \sum_{n \in \ZZ} e^{2 \pi ia (n+\frac{c}{2} )}
f_\epsilon((n+c)x) \right)
x^{s-1} dx & \mbox{if} & \ep = (-1)^k, \\ [+.2in]0 & \mbox{if} & \ep
= (-1)^{k+1},
\end{array}
\right.
 \eeq 
 by the same proof as \eqn{JL212} above.
Thus we have $\Psi_k (f;s,a,c) = \Psi_k(f_\ep;s,a,c)$
with $\eps~=~(-1)^k$. The Poisson summation
formula applied to Schwartz functions gives, for $x \in \RR_{>0}$
and any $a,c \in \RR$,
 \beql{JL218} 
 \sum_{n \in \ZZ} e^{2 \pi ian} f((n+c)x) = 
 \frac{1}{|x|} \sum_{m \in \ZZ} e^{2 \pi ic (m-a)} \sF f
\left( \frac{m-a}{x} \right)~. 
\eeq 
Substituting this into the
definition of $\Psi_k (s,a,c;f)$ yields, for $\ep = (-1)^k$, that
$$
\Psi_k (f; s,a,c)= \Psi_k (f_\ep; s,a,c) = 2 \int_0^1 \left( e^{\pi iac} \sum_{m \in
\ZZ} e^{2 \pi ic (m-a)} \sF f_{\epsilon} \left( \frac{m-a}{x}
\right)\right) x^{s-1} dx\,.
$$
Making the change of variable $y= 1/x$ gives rise to
 \beql{JL219}
\Psi_k (f_\ep; s,a,c) = 2 \int_1^\infty \left( \sum_{m \in \ZZ} e^{2
\pi ic (m - \frac{a}{2} )} \sF f_{\epsilon} ((m-a) y)\right) y^{-s}
dy \,. \eeq Letting $m \to -m$ yields, since $\ep = (-1)^k$,
\begin{eqnarray}\label{JL220}
\Psi_k ( f_\ep;  s,a,c) & = &
2 \int_1^\infty \left( \sum_{m \in \ZZ} e^{2 \pi ic (-m- \frac{a}{2} )}
\sF f_\ep
((-m-a) y)\right) y^{-s} dy \nonumber \\
& = & (-1)^k e^{ - \pi ia}
\left( 2 \int_1^\infty \sum_{m \in \ZZ} e^{2 \pi i (1-c) (m+ \frac{a}{2} )}
\sF f_\ep ((m+a) y) \right) y^{-s} dy \nonumber \\
& = & (-1)^k e^{- \pi ia} \Phi_k ( \sF f_\ep;   1-s, 1-c, a) ~.
\end{eqnarray}
Since $\Phi_k (f; s,a,c) = \Phi_k (f_{\ep}; s,a,c)$, by a similar proof
to \eqn{JL212},
we obtain from \eqn{JL220} that
\beql{JL220a}
\Psi_k (f; s,a,c) = (-1)^k e^{- \pi ia} \Phi_k
( \sF f;  1-s , 1-c, a)
\eeq
for all $f \in \sS(\RR)$.
Substituting this in \eqn{JL215} yields the desired result.~~~$\bsq$ \\

We now deduce Theorem ~\ref{th11} using  Lemma~\ref{le21} 
with specific test functions.\\


\paragraph{\bf Proof of Theorem~\ref{th11}.}
The functional equations are associated to test functions which
are self-reciprocal under the Fourier transform, in the
sense that $\sF f(x) = c f(x)$, for some constant $c$, necessarily a fourth
root of unity.
We  consider  the test functions
$$
\phi_0(x) := e^{- \pi x^2} \in \sS^+  ~~\mbox{and}~~
\phi_1(x) := x~e^{- \pi x^2} \in \sS^-.
$$
These satisfy the self-reciprocal 
Fourier-transform formula $ \sF\phi_0(x)=e^{- \pi x^2}= \phi_0(x)$ and
$ \sF \phi_1(x)=-i\phi_1(x)$.

Suppose that $0 < a < 1$, $0 < c < 1$ and $\Re (s) > 1$. Since
$\sM_0 (e^{- \pi x^2})(s) = \pi^{- s/2} \Gamma (s/2)$, the 
definition of $F_k (f;s,a,c)$ gives
\beql{371}
\hat{L}^{+}(s,a,c) = \pi^{-s/2} \Gamma (s/2) (\zeta (s,a,c) +
e^{- 2 \pi ia} \zeta (s, 1-a, 1-c ))
= e^{- \pi iac} F_0 (\phi_0; s,a,c)   ~.
\eeq
Similarly since
$\sM_1 ( xe^{- \pi x^2})(s) = \pi^{- \frac {s+1}{2}}
\Gamma \left( \frac{s+1}{2} \right)$,
one finds that
\beql{372}
\hat{L}^{-}( s,a,c)=  \pi^{- \frac{s+1}{2}}
\Gamma \left( \frac{s+1}{2} \right) (\zeta (s,a,c) - e^{- 2 \pi
ia} \zeta (s,1-a, 1-c)) = e^{- \pi iac} F_1 (\phi_1; s,a,c) ~.
\eeq
 Lemma ~\ref{le21} gives the analytic continuation of 
$F_k(f; s , a , c)$ to an entire function of $s$, which shows for 
$0< a<1, 0< c<1$
 that each of 
$\hat{L}^{\pm}(s, a, c)$ extends to an  entire function of $s$.

The functional equations  follow from
Lemma~\ref{le21} using  the definition of $\Phi_k(f; s, a, c)$ and the Fourier transform
formulas.  Namely
\begin{eqnarray*}
\hat{L}^{+}( s, a ,c) & =& e^{-\pi i a c} F_0(\phi_0; s , a, c) \\
&=& e^{-\pi i a c} \left(  \Phi_0(\phi_0; s, a, c) + e^{-\pi i a} \Phi_0(\phi_0; 1-s, 1-c, a)
\right) \\
&=& e^{-\pi i ac} \left( \Phi(\phi_0; s, a, c) + e^{-\pi i(a+1-c)} \Phi(\phi_0; s, 1-a, 1-c) \right) \\
&& + e^{-\pi i  a  c} \left( e^{-\pi i a} \Phi(\phi_0; 1-s, 1-c, a)+ e^{-\pi i c}\Phi(\phi_0; 1-s, c, 1-a) 
\right) ,
\end{eqnarray*}
which agrees with 
\begin{eqnarray*}
e^{-2 \pi i  ac} &&\hat{L}^{+}(1-s, 1-c, a)\\
 &=& e^{-2 \pi i  ac} e^{- \pi i (1-c)a} F_{0} (\phi_0; 1-s, 1-c, a) \\
&=& e^{- \pi i  ac} e^{- \pi i a} \left(
 \Phi_0(\phi_0; 1-s, 1-c, a) + e^{-\pi i (1-c)} \Phi_0(\phi_0; s, 1-a, 1-c) \right)\\
&=& e^{-\pi i a c}  e^{- \pi i a}  \left(   \Phi( \phi_0; 1-s, 1-c, a) + 
e^{-\pi i (-c-a)} \Phi(\phi_0; 1-s, c, 1-a) \right)\\
&&
+  e^{-\pi i a c- \pi i a} e^{-\pi i (1-c)} \left( \Phi( \phi_0; s, 1-a , 1-c) 
+ e^{-\pi i (-a + c -1)} \Phi(\phi_0; s, a ,c)\right)  \\
&=& e^{- \pi i a c} \left(
 e^{- \pi i a}  \Phi( \phi_0; 1-s, 1-c, a) + e^{-\pi i c} \Phi(\phi_0;1-s, c, 1-a) \right) \\
 &&  + e^{- \pi i  a c} \left( e^{-\pi i (a+1-c)}  \Phi( \phi_0; s, 1-a , 1-c)+ 
  \Phi(\phi_0; s, a ,c) \right).
\end{eqnarray*}
A similar derivation holds for $\hat{L}^{-}(s, a,c) = i e^{-2\pi i a c}\hat{L}(1-s, 1-c, a),$
using $\phi_1(x)$.
~~~$\bsq$ \\

\paragraph{\bf Remark.}   
Weil's proof \cite[p. 57]{We76} of the four term functional equations 
in Theorem~\ref{th11} uses
Kronecker's summation formula, which is a variant of the 
Poisson summation formula, applied directly to 
$\phi_0$ and $\phi_1$.  The use,
explicitly or implicitly,  
of the Poisson summation 
seems inevitable, in that
 it is known to be 
equivalent to the functional equation for $\zeta(s)$, see
Patterson \cite[Sec. 2.14]{Pa88}.

%

\section{Extended zeta integrals}
\setcounter{equation}{0}

We now  extend the  ``zeta integrals" to apply to 
all $(a,c) \in \RR \times \RR$. Using the formulation of Lemma~\ref{le21} for
the ``zeta integral" $F_k(f; s, a, c)$ that
$$
 F_k (f; s,a,c) =
\Phi_k (f; s,a,c) + (-1)^k e^{- \pi ia} \Phi_k ( \sF f;  1-s, 1-c,a)
$$
it suffices to extend the definition of  $\Phi_k (f; s,a,c)$ to apply to all
$(a,c) \in \RR \times \RR$.
Now the  function $\Phi (f; s,a,c)$ given by \eqn{JL213b} is well-defined
for all
$(a,c) \in \RR \times \RR$, provided $\Re(s) < 0$.
When $c \not\in \ZZ$, the integral representation \eqn{JL213b}
converges for all $s \in \CC$ and
defines $\Phi (f; s,a,c)$ as an entire
function of $s$, and the same property is inherited by $\Phi_k (f; s,a,c)$,
for $k = 0,1.$
In the remaining case
 $c \in \ZZ$ and $\Re (s) < 0$, the term $n = -c$ in \eqn{JL213b}
contributes 
\beql{JL221}
 e^{-\pi iac} f(0) \int_1^\infty x^{s-1}
dx = -e^{-\pi iac} f(0) \frac{1}{s} , 
\eeq
while the remaining terms define an entire function of $s$. The
functions $\Phi_k (f; s,a,c)$ in \eqn{JL213a} have two such terms,
the second with weight $(-1)^k$. Since the right side of
\eqn{JL221} analytically continues to  $s \in \CC$, this motivates
the following definition.
\begin{defi}\label{de21}
{\rm
For  $k=0$ or $1$, $(a,c) \in \RR \times \RR$ and $s \in \CC$, define
for each $f(x) \in \sS(\RR)$,
\begin{eqnarray}\label{JL222}
\lefteqn{\Phi_k (f; s,a,c) := \int_1^\infty \left(
\sum_{n \in \ZZ \atop n \neq -c} e^{2 \pi ia (n+ \frac{c}{2} )}
f((n+c) x) \right) x^{s-1} dx}~~~~~~~~~~~~~~ \nonumber \\
&& \quad + (-1)^k e^{- \pi i (a+1-c)} \int_1^\infty \left(
\sum_{m \in \ZZ \atop m\neq c-1} e^{2 \pi i (1-a) (m+ \frac{1-c}{2} )}
f((m+1-c) x) \right) x^{s-1} dx \nonumber \\
&&\quad - \delta_{\ZZ}(c)~ e^{-\pi iac} (1+ (-1)^k) f(0) \frac{1}{s} ~,
\end{eqnarray}
in which
\beql{JL223}
\delta_{\ZZ}(c) =1 \quad\mbox{if}\quad c \in \ZZ \quad\mbox{and~ equals}
\quad 0 \quad\mbox{otherwise} ~.
\eeq
}\end{defi}

For fixed $a$ and $c$, the function $\Phi (f;s,a,c)$
 is a meromorphic function
in $s$. It is holomorphic except for a possible simple pole at
$s=1$, which occurs if and only if $k=0$, $c \in \ZZ$  and $f(0)
\neq 0$. We make the following definition, which in view of
Lemma~\ref{le21} extends Definition~\ref{de20} to $(a,c) \in \RR
\times \RR,$ and to all $s \in \CC.$
\begin{defi}\label{de22}
{\rm Let $(a,c) \in \RR \times \RR$ and $f \in \sS(\RR )$. Then
for $k~=~0,~1$ and all $s \in \CC$, set }
\beql{JL225}
 F_k (f;s,a,c) := \Phi_k (f; s,a,c) + (-1)^k e^{-\pi i a} \Phi_k(\sF f; 1-s, 1-c, a). 
 \eeq 
 \end{defi}

We note that the dependence of $F_k(f; s,a,c)$ on the test function
$f$ is still entirely in the Mellin transform factor $\sM_k(f)(s)$,
as in Lemma~\ref{nle21}.

\begin{theorem}\label{JLth21} {\rm (Lerch Zeta Integral Functional Equations)}
Let  $(a,c) \in \RR \times \RR$ and $f(x) \in {\mathcal S}( \RR )$.
Then the functions
$$
F_k (f; s,a,c) := \Phi_k (f; s,a,c) + (-1)^k e^{- \pi ia} \Phi_k ( \sF f; 1-s, 1-c, a)
$$
are meromorphic functions of $s$, having the following
properties. 

\begin{itemize}
\item[(i)] The function $F_1(f; s, a ,c)$ is an entire function of $s$ for
 all $(a, c) \in \RR \times \RR$. The function $F_0 (f; s,a,c)$
is an entire function of $s$ for  $(a, c) \in \RR \times \RR$ with
both $a, c$ non-integers. Otherwise it is meromorphic, 
with possible  simple poles at $s=0, 1$. A pole may occur 
at $s=0$ only if $a$ is an integer, and at $s=1$ only if $c$
is an integer, with residues  given by 
\begin{eqnarray}
Res_{s=0} F_0 (f; s,a,c) & = & -2e^{-\pi i a c} f(0)
~~~~~\mbox{~~~~if}~~ c \in \ZZ, \label{JP223a}\\
Res_{s=1} F_0 (f; s,a,c) & = & ~2 e^{\pi i a c} \sF f(0)
~~~~~\mbox{~~~~~if}~~ a\in \ZZ.\label{JP223b}
\end{eqnarray}

\item[(ii)]
The functions $F_k (f; s,a,c)$
are twisted periodic in the $a$-variable and
$c$-variable, satisfying
\begin{eqnarray}
\label{JP224}
F_k (f; s, a+1 , c) & = & e^{\pi ic} F_k (f; s,a,c), \\
\label{JP225} F_k (f; s,a,c+1) & = & e^{- \pi ia} F_k (f; s,a,c).
\end{eqnarray}

\item[(iii)]
The functions  $F_k (f; s,a,c)$ satisfy the functional equations
\beql{JP226}
F_k (f; s,a,c) = (-1)^k e^{- \pi ia} F_k ( \sF f ;   1-s, 1-c, a) ~.
\eeq
\end{itemize}
\end{theorem}

\paragraph{\bf Proof.}
(i) For $(a,c) \in \RR \times \RR$ the definition \eqn{JL222} gives
\begin{eqnarray*}
F_k (f; s,a,c)  &= & \delta_{\ZZ}(a)~ e^{\pi iac} (1+ (-1)^k )
\sF f(0) \frac{1}{s-1}\\
&&~~~ - \delta_{\ZZ}(c)~ e^{- \pi iac} (1+ (-1)^k) f(0)
\frac{1}{s} + R_k (f; s,a,c) ~,
\end{eqnarray*}
in which $R_k (f; s,a,c)$ is an entire function of $s$.
This establishes the analytic properties of
$F_k (f; s,a,c)$ in the $s$ variable.

(ii). For $\ep = (-1)^k$ we have
\begin{eqnarray}\label{JP227}
\Phi_k (f_\ep ; s,a+1,c) & = &
2 \int_1^\infty \left( \sum_{n \in \ZZ} e^{2 \pi i (a+1) (n+ \frac{c}{2} )}
f_\ep ((n+c)x) \right) x^{s-1} dx \nonumber \\
& = & e^{\pi ic} \Phi_k (f_\ep ; s,a,c) ~,
\end{eqnarray}
and, similarly,
\begin{eqnarray}\label{JP228}
\Phi_k ( f_\ep ;  s,a,c+1) & = & 2 \int_1^\infty
\left( \sum_{n \in \ZZ} e^{2 \pi ia(n+ \frac{c+1}{2} )}
f_\ep ((n+c+1)x) \right) x^{s-1} dx \nonumber \\
& = & e^{- \pi ia} \Phi_k (f_\ep ; s,a,c) ~.
\end{eqnarray}
These identities give, using \eqn{JL210a},
\begin{eqnarray*}
F_k (f; s,a+1,c) & = & F_k (f_\ep; s, a+1, c) \\
& = & \Phi_k (f_\ep ; s,a+1, c) + (-1)^k
e^{- \pi i(a + 1)} \Phi (\sF f_\ep; 1-s, 1-c, a+1) \\
& = & e^{\pi ic} \Phi (f_\ep ; s,a,c) + (-1)^k e^{- \pi i(a+1)}
e^{- \pi i (1-c)} \Phi (\sF f_\ep; 1-s, 1-c, a)
\\
& = & e^{ \pi ic} F_k ( f_\ep; s,a,c) \\
& = & e^{\pi ic} F_k (f; s,a,c) ~.
\end{eqnarray*}
In a similar fashion they also give
$$F_k (f; s,a,c+1) = e^{- \pi ia} F_k (f; s,a,c). $$

 (iii) We first show that
\beql{JP229} 
\Phi_k (f; s,1-a,1-c) = (-1)^k
e^{\pi i (a+1-c)} \Phi_k (f; s,a,c) ~.
 \eeq

We have, for $\ep = (-1)^k$, that
\begin{eqnarray*}
\Phi_k (f; s,1-a, 1-c) & = & \Phi (f_{\ep}; s,1-a, 1-c) \\
& = & 2
\int_1^\infty \left( \sum_{n \in \ZZ} e^{2 \pi i (1-a) (n+
\frac{1-c}{2} )} f_\ep ((n+1-c)x) \right) x^{s-1} dx \,.
\end{eqnarray*}
Replacing $n$ with $-n-1$ yields
\begin{eqnarray*}
\Phi_k (f_{\ep};s,1-a,1-c) & = & 2 \int_1^\infty \left( \sum_{n \in
\ZZ} e^{2 \pi i (1-a) (-n-(\frac{1+c}{2} ))}
f_\ep ((-n-c)x) \right) x^{s-1} dx \\
& = & 2 (-1)^k e^{\pi i (a+1-c)} \int_1^\infty
\sum_{n \in \ZZ} e^{2 \pi ia (n+\frac{c}{2})}
f_\ep ((n+c)x) x^{s-1} dx \\
& = & (-1)^k e^{\pi i (a+1-c)} \Phi_k (f_{\ep}; s,a,c) = (-1)^k
e^{\pi i (a+1-c)} \Phi_k (f; s,a,c) ~,
\end{eqnarray*}
which is \eqn{JP229}.

We deduce the functional equation \eqn{JP226} as follows.
For $\ep = (-1)^k$, we have
\begin{eqnarray}\label{JP230}
F_k (\sF f; 1-s, 1-c, a)
& =  & F_k(\sF f_\ep; 1-s , 1-c, a) \nonumber  \\
& = & \Phi_k (\sF f_\epsilon; 1-s, 1-c, a ) +
(-1)^k e^{- \pi i(1-c)} \Phi_k (\sF \sF f_\ep ; s, 1-a, 1-c ) \nonumber \\
& = &(-1)^k e^{\pi ia}
((-1)^k e^{- \pi ia} \Phi_k (\sF f_\ep; 1 -s, 1-c, a) \nonumber \\
&&~~~~~~+ \ep e^{- \pi i(a+1-c)} \Phi_k (f_\ep;  s, 1-a, 1-c) ) \nonumber \\
& = & (-1)^k e^{\pi ia} ((-1)^k e^{- \pi ia} \Phi_k (\sF f_\ep; 1-s, 1-c, a)
+ \Phi_k (f_\ep; s,a,c)) \nonumber \\
& = &(-1)^k e^{\pi ia} F_k (f_{\ep}; s,a,c) \nonumber \\
& = & (-1)^k e^{\pi ia} F_k (f; s,a,c) ~, \nonumber
\end{eqnarray}
in which we used $\sF \sF f_\ep (x) = f_\ep (-x) = \ep f_\ep (x)$, and
\eqn{JP229} was used at the third step.~~~$\bsq$ \\

%
%
%
%
%

\section{Extended Lerch zeta function}
\setcounter{equation}{0}

We next enlarge the definition of the functions $L^{\pm}( s ,a,c)$ to
$(a,c) \in \RR \times \RR$, using certain extended zeta integrals from \S4.  
Then we use  these functions to define an extended Lerch zeta function
$\zeta_{*}(s, a,c)$ to the same domain. We emphasize that this extended
function is not obtainable by analytic continuation in the $a$- anc $c$-variables. 
Indeed our  results  show that  these extended functions necessarily
have discontinuities in the $(a,c)$ variables at integer values of $a$ and $c$,
for some ranges of $s$; their virtue is that they preserve the twisted periodicity
relations given in Theorem \ref{th12}.

Using the test functions $\phi_0(x)= e^{- \pi x^2}$ and $\phi_1(x)= x e^{- \pi x^2}$
we obtained in \S3 for $(a,c) \in \Box^{\circ}$ the relations 
$$
\hat{L}^{\pm}(s, a, c) = e^{-\pi i a c} F_k( \phi_k; s , a ,c),    
$$
where $k=0, 1$ and $\pm= (-1)^k$, in \eqn{371}, \eqn{372}. 
The results of \S4 now define the right side for all $(a,c) \in \RR \times \RR$, 
which motivates the following definition.

\begin{defi}\label{de51a}
{\em For  all $(a, c) \in \RR \times \RR$ and all $s \in \CC$ the
 ({\em extended}) {\em completed Lerch zeta functions} $L_{\ast}^{\pm}( s,a,c)$
are given by 
\beql{JP233}
\hat{L}_{\ast}^{+}(s,a,c) := e^{- \pi iac} F_0 (\phi_0; s,a,c),
\eeq 
\beql{JP234}
 \hat{L}_{\ast}^{-}(s,a,c) := e^{- \pi iac} F_1(\phi_1; s,a,c) ~, 
\eeq 
where the right hand sides are meromorphic functions of $s \in \CC$,
using  Theorem~\ref{JLth21}. }
\end{defi}

Using this definition we reverse-engineer  an extension
of the Lerch zeta function  $\zeta (s,a,c)$ on
the domain $\Box^{\circ}$ to a function $\zeta_{\ast}(s,a,c)$ on
$(a,c) \in \RR\times \RR$.
To do this, view  the equation \eqn{109b}
as a pair of linear equations in $\zeta (s,a,c)$
and $\zeta (s,1-a,1-c)$, and eliminate
$\zeta (s,1-a,1-c)$ to obtain
\beql{JP235}
\zeta (s,a,c) = \frac{1}{2} \left(
\frac{\pi^{s/2}}{\Gamma (\frac{s}{2})} \hat{L}^{+}( s,a,c)  +
\frac{\pi^{(s+1)/2}}{\Gamma (\frac{s+1}{2} )}
\hat{L}^{-}( s,a,c) \right),
\eeq
which is valid for $0 < a < 1$, $0 < c < 1$, and $\Re (s) > 1$.
The right side now extends by Definition \ref{de51a}  to  all $(a, c) \in \RR \times \RR,$
which motivates to the following definition.


\begin{defi}\label{Jde24}
 {\rm For $(a,c) \in \RR \times \RR$ and $s \in \CC$, the
{\em extended Lerch zeta function} $\zeta_\ast (s,a,c)$ is defined
by 
\beql{JP236}
 \zeta_\ast (s,a,c) := \frac{1}{2} \left( \frac{\pi^{s/2}}{\Gamma ( \frac{s}{2})}
\hat{L}_{\ast}^{+}(s,a,c) +
 \frac{\pi^{(s+1)/2}}{\Gamma ( \frac{s+1}{2})} \hat{L}_{\ast}^{-}( s,a,c) \right) ~. 
\eeq }
\end{defi}

The function $\zeta_\ast (s,a,c)$ is clearly meromorphic in $s$,
and it is holomorphic in $\CC$ except for a possible simple
pole at $s=1$ coming from $ \hat{L}^{+}( s,a,c)$. To see this,
observe that any pole at $s=0$ of $ \hat{L}^{+}( s,a,c)$
is cancelled by the
zero of $\frac{1}{\Gamma (\frac{s}{2})}$
at $s=0$, and the remaining terms on the right side of
\eqn{JP236} are entire functions.
 Also \eqn{JP235} shows that for
$0 < a < 1$ and $0 < c < 1$ we have 
\beql{JP236a}
 \zeta_\ast
(s,a,c) = \zeta (s,a,c) ~.
 \eeq 
The following  lemma establishes the Dirichlet
series formula \eqn{sp111b} for $\zeta_\ast (s,a,c)$, as well as
formulae for $\hat{L}^{\pm}(s,a,c)$.
\begin{lemma}\label{Jle22}
For $\Re (s) > 1$  and all $(a,c) \in \RR \times \RR$, there holds
\begin{eqnarray}~\label{JP237}
 \hat{L}_{\ast}^{+}( s,a,c) & = & \pi^{- \frac{s}{2}} \Gamma \left(
\frac{s}{2} \right) \left( \sum_{n \in \ZZ \atop n \neq -c} e^{2 \pi
ina} |n+c|^{-s} \right) ~, \\
\hat{L}_{\ast}^{-}(s,a,c) &= & \pi^{- \frac{s+1}{2}} \Gamma \left(
\frac{s+1}{2} \right) \left( \sum_{n \in \ZZ \atop n \neq -c} sgn
(n+c) e^{2 \pi ina} |n+c|^{-s} \right) ~.~\label{JP238}
\end{eqnarray}
In consequence, on the same domain,
\beql{JP239}
\zeta_\ast (s,a,c) = \sum_{n+c > 0} e^{2 \pi ina} (n+c)^{-s} ~.
\eeq
\end{lemma}

\pf
For $s \in \CC$, evaluating $ \hat{L}_{\ast}^{+}( s,a,c)$ from its
definition \eqn{JP233}, using \eqn{JL225} and \eqn{JL222},
yields
\begin{eqnarray}\label{JP240}
\frac {1}{2} \hat{L}_{\ast}^{+}( s,a,c) & = & \int_1^\infty
\left( \sum_{n \in \ZZ \atop n \neq -c}
e^{2 \pi ian} e^{- \pi (n+c)^2 x^2} \right) x^{s-1} dx \nonumber \\
&& ~~+ \int_1^\infty
\left( \sum_{n \in \ZZ \atop n \neq a}
e^{2 \pi ic (n-a)} e^{- \pi (n-a)^2 y^2} \right) y^{-s} dy \nonumber \\
&&~~+ \delta_{\ZZ}(a) \frac{1}{s-1}
 - \delta_{\ZZ}(c)~e^{-2\pi i a c} \frac{1}{s} ~.
\end{eqnarray}
Since $\Re (s) > 1$ the $\delta_{\ZZ}(a)$-term
may be absorbed as $n=a$ in the second integral. Applying Poisson
summation \eqn{JL218} to the second integral and letting 
$x =\frac {1}{y}$ yields, for $\Re (s) > 1$,
\begin{eqnarray}
 \int_1^\infty \left( \sum_{n \in \ZZ} e^{2 \pi ic(n-a)}
e^{- \pi (n-a)^2 y^2} \right) y^{-s} dy
 & = & \int_0^1 \left( \sum_{n \in \ZZ \atop n \neq -c}
e^{2 \pi ian} e^{-\pi (n+c)^2 x^2} \right) x^{s-1} dx \nonumber \\
&& +~~ \delta_{\ZZ}(c) e^{-2 \pi i a c}\int_0^1 y^{s-1} dy ~.
\end{eqnarray}
The last integral is $\frac{1}{s}$ when $\Re (s) > 0$, and substituting
this in \eqn{JP240} yields
$$
\frac {1}{2}\hat{L}_{\ast}^{+}(s,a,c) = \int_0^\infty
\left( \sum_{n \in \ZZ \atop  n \neq -c}
e^{2 \pi ian} e^{- \pi (n+c)^2 x^2} \right) x^{s-1} dx ~,
$$
in which  the pole terms at $s= 0$ when $c \in \ZZ$ cancel out.
Since $\Re (s) > 1$ we can interchange summation with the integral,
and using the change of
variables $y= | n+c| x$ yields \eqn{JP237}.

We can perform a similar sequence of steps on $\hat{L}_{\ast}^{-}(s,a,c)$
starting from
\begin{eqnarray}\label{JP241}
\frac {1}{2} \hat{L}_{\ast}^{-}( s,a,c) & = & \int_1^\infty
\left( \sum_{n \in \ZZ}
e^{2 \pi ian} (n + c)e^{- \pi (n+c)^2 x^2} \right) x^{s} dx \nonumber \\
&& ~~- i \int_1^\infty
\left( \sum_{n \in \ZZ}
e^{2 \pi ic (n-a)}(n - a) e^{- \pi (n-a)^2 y^2} \right) y^{1-s} dy,
\end{eqnarray}
which is valid for $s \in \CC$, to obtain \eqn{JP238}.
Substituting these two formulas in the definition
\eqn{JP236} of $\zeta_\ast (s,a,c)$ yields \eqn{JP239}.~~~$\bsq$\\

We now deduce Theorem~\ref{th12} from Theorem~\ref{JLth21}
with special  test functions. \\

\paragraph{\bf Proof of Theorem~\ref{th12}.}
Using Lemma \ref{Jle22} we see that for $\Re(s)>1$ the
definition \ref{Jde24} for $\zeta_\ast (s,a,c)$
 agrees with the formula \eqn{sp111b},
\[
\zeta_{\ast}(s, a,c) = \sum_{n+c >0} e^{2\pi in a}(n+c)^{-s}.
\]
For $\Re(s)>1$ this  function is manifestly
discontinuous in the $c$-variable, it jumps
when $c$ is an integer. In particular, when $\Re(s)>1$, we observe that
\[
\zeta_{\ast}(s, a,c+1) = \sum_{n+c+1 >0} e^{2\pi in a}(n+c+1)^{-s}= e^{-2\pi i a}
\sum_{m+c>0} e^{2\pi i m a}(m+c)^{-s},
\]
which is a special case of property (ii).

Properties (i) - (iii) of Theorem \ref{th12} follow using
definition \ref{Jde24}  by simple
calculations from properties (i) and (ii) of
Theorem \ref{JLth21}, applied to the  test functions
$\phi_0(x) = e^{- \pi x^2}$ and $\phi_1(x) = x e^{- \pi x^2}$.
The changes in the form of the twisted periodicity (ii) and the functional
equation (iii) compared to those in Theorem ~\ref{JLth21}
 come from the $e^{- \pi iac}$ factor in
$\hat{L}_{\ast}^{\pm}( s,a,c) = e^{ - \pi iac} F_k (\phi_k; s,a,c)$, for
$k= 0, 1$, with $\pm = (-1)^k$.~~~$\bsq$ \\

We conclude this section by deducing Lerch's transformation formula
for the extended Lerch zeta function $\zeta_{\ast} (1-s,a,c)$
from Theorem \ref{th12}. 
\begin{theorem}\label{th51}
{\rm (Extended Lerch's transformation formula)} For $(a, c) \in \RR \times
\RR$ and all $s \in \CC$,
\beql{sp34}
 \zeta_{\ast} (1-s,a,c) =  (2 \pi )^{-s} \Gamma (s) 
 \left\{ e^{\frac{\pi is}{2}} e^{- 2 \pi iac} \zeta_{\ast} (s,1-c,a) + e^{-
\frac{\pi is}{2}} e^{ 2 \pi ic (1-a)} \zeta_{\ast} (s,c,1-a) \right\}.
\eeq
\end{theorem}

\pf
We combine the functional equations of $\hat L^{\pm}( s, a,c)$
to eliminate one term. We obtain
\begin{eqnarray}\label{sp34a}
2~ \zeta_{\ast}(s, a, c)  =  e^{-2 \pi i a c} \{
( \frac {\pi^{-\frac {1-s}{2}} \Gamma(\frac {1-s}{2})}{\pi^{-\frac {s}{2}}
\Gamma(\frac {s}{2})} +
i~ \frac {\pi^{-\frac {2-s}{2}}
\Gamma(\frac {2-s}{2})}{\pi^{-\frac {1 +s}{2}}
\Gamma(\frac {1+s}{2})} ) \zeta_*(1-s, 1-c,a) \nonumber\\
~~~~~+e^{2 \pi i c}
( \frac {\pi^{-\frac {1-s}{2}} \Gamma(\frac {1-s}{2})}{\pi^{-\frac {s}{2}}
\Gamma(\frac {s}{2})} -
i~ \frac {\pi^{-\frac {2-s}{2}}
\Gamma(\frac {2-s}{2})}{\pi^{-\frac {1 +s}{2}}
\Gamma(\frac {1+s}{2})} ) \zeta_*(1-s, c, 1- a)\}.
\end{eqnarray}
Then, replacing $s$ by $1 - s$ throughout and judiciously applying
the two gamma function identities 
\beql{sp35} 
\Gamma (s) \Gamma(1-s) = \frac{\pi}{\sin \pi s} 
\eeq
 and 
 \beql{sp36} 
 \Gamma \left(\frac{s}{2} \right) \Gamma \left( \frac{1+s}{2} \right)=
 \sqrt{2 \pi} 2^{\frac {1} {2} -s} \Gamma (s)
\eeq
 eventually yields \eqn{sp34}.~~~$\bsq$

%
%
%
%
%

\section{Boundary cases of functional equation}
\setcounter{equation}{0}

In this section we study the behavior of $\zeta_\ast (s,a,c)$ as $(a,c)$
moves from the interior of the unit square to the boundary, for fixed $s$;
these formulas extract terms producing   discontinuities at these boundaries, 
whose nature depends on
the value of $\Re(s)$.
\begin{theorem}\label{th31}
{\em (Lerch Boundary Limits)}
 For fixed $s \in \CC \smallsetminus  \ZZ$
and $(a,c) \in \Box^\circ$, the following hold.
\begin{itemize}
\item[(i)]
For $0 < a' < 1$ and $(a,c) \to (a',1^-)$, \beql{P31} \zeta_\ast
(s,a,c)~~ \to ~~\zeta_\ast (s, a', 1). \eeq
\item[(ii)]
For $0 < a' < 1$ and $(a,c) \to (a',0^+)$, \beql{P32} \zeta_\ast
(s,a,c) - c^{-s}~~ \to~~ \zeta_\ast (s, a', 0). \eeq
\item[(iii)]
For $0 < c' < 1$ and $(a,c) \to (1^-, c')$,
\begin{eqnarray}\label{P33}
 \zeta_\ast (s,a,c) &- &\frac{1}{2} \frac{\pi^{\frac{s-1}{2}}
\Gamma \left( \frac{1-s}{2} \right)}{\pi^{- \frac{s}{2}}
\Gamma \left( \frac{s}{2} \right)}
\left( 1- i \cot \frac{\pi s}{2} \right) e^{2 \pi i(1-a)c} (1-a)^{s-1}\\
~~ && \mbox{~~~~~~~~~~~~}~~~~~~~~~~~~~~~~~~\to~~ \zeta_\ast (s,1,c'). \nonumber
\end{eqnarray}
\item[(iv)]
For $0 < c' < 1$ and $(a,c) \to (0^+,c')$,
\begin{eqnarray}\label{P34}
\zeta_\ast (s,a,c) - \frac{1}{2} \frac{\pi^{\frac{s-1}{2}}
\Gamma \left( \frac{1-s}{2} \right)}{\pi^{- \frac{s}{2}}
\Gamma \left( \frac{s}{2} \right)}
\left( 1 + i \cot \frac{\pi s}{2} \right) e^{-2 \pi iac} a^{s-1}
~~ \to ~~\zeta_\ast (s,0, c').
\end{eqnarray}
\end{itemize}
Furthermore (i) and (ii) are valid for integer $s \le 0$, while
(iii) and (iv) are valid for integer $s \ge 1.$
\end{theorem}

To prove this result we begin with a preliminary lemma.
\begin{lemma}\label{le31}
For $0 < t < \infty$ and $s \in \CC$,
\beql{P35}
\int_1^\infty e^{- \pi t^2 x^2} x^{s-1} dx = - \frac{1}{s} e^{- \pi
t^2} + \pi^{-s/2} \Gamma \left( \frac{s}{2} \right) t^{-s} + h(t,s)
\,, \eeq where $h(t,s)$ is continuous in $t$ and meromorphic for $s
\in \CC$, with poles only at $s \in \{0, -2, -4, ...\}$, and for $s
\not\in 2\ZZ_{\le 0}$, 
\beql{P36}
 \lim_{t \to 0^+} h(t,s) =0.
  \eeq
Under the same conditions 
\beql{P37} 
\int_1^\infty te^{-\pi t^2 x^2}x^{s} dx = 
\pi^{- \frac{s+1}{2}} \Gamma \left( \frac{s+1}{2} \right) t^{-s} + k(t,s), 
\eeq 
where $k(t,s)$ is continuous in $t$ and
meromorphic in $s \in \CC$, with poles only at 
$s \in \{-1, -3, -5,...\}$, and for $s \notin -1+2\ZZ_{\le 0}$,
 \beql{P38}
  \lim_{t \to 0^+} k(t,s) =0. 
  \eeq
\end{lemma}

\pf
It suffices to prove the statement in \eqn{P36} for
 $s \in \CC \smallsetminus 2\ZZ_{\le 0}$
 since both sides of \eqn{P35} are
meromorphic in $s$. Let $m \in \ZZ_{\ge 0}$ be such that $\Re (s+2m)
> 0$. Integrating by parts repeatedly gives
\begin{eqnarray*}
 \int_1^\infty e^{- \pi t^2 x^2} x^{s-1} dx & =  &
e^{- \pi t^2 x^2} \frac{x^s}{s}
\Bigl|_1^\infty +
\frac{2 \pi t^2}{s} \int_1^\infty e^{- \pi t^2 x^2} x^{s+1} dx \\
& = & - \frac{e^{- \pi t^2}}{s} +
\frac{2 \pi t^2}{s} \int_1^\infty e^{- \pi t^2 x^2} x^{s+1} dx \\
& = & - \sum_{k=0}^m \frac{(2 \pi)^k t^{2k}
e^{- \pi t^2}}{s(s+2) \cdots (s+2k)} \\
&& ~~~~~~ + 
\frac{(2 \pi )^{m+1} t^{2m+2}}{s(s+2) \cdots (s+2m)} \int_1^\infty
e^{- \pi t^2 x^2} x^{s+2m+1} dx \,.
\end{eqnarray*}
Now write the last integral as $\int_0^\infty - \int_0^1$ and note
that $\int_0^1 e^{- \pi t^2 x^2} x^{s+2m+1} dx$ converges absolutely
since $\Re (s+2m +1) > 0$. We obtain \beql{P39} \int_1^\infty e^{-
\pi t^2 x^2} x^{s-1} dx = - \frac{1}{s} e^{- \pi t^2} + h(t,s) +
\frac{(2 \pi )^{m+1} t^{2m+2}}{s(s+2) \cdots (s+2m)} \int_0^\infty
e^{-xt^2 x^2} x^{s+2m+1} dx \,, \eeq where we define \beql{P310}
h(t,s) := - \sum_{j=1}^m \frac{(2 \pi )^j t^{2j} e^{- \pi
t^2}}{s(s+2) \cdots (s+2j)} ~~-~~ \frac{(2 \pi )^{m+1}
t^{2m+2}}{s(s+2) \cdots (s+2m)} \int_0^1 e^{- \pi t^2 x^2} x^{s+
2m+1} dx \,.
\eeq

Using the change of variable $y = tx$, we can rewrite the last term in
\eqn{P39} as
\begin{eqnarray}
\label{P311} \frac{(2 \pi )^{m+1} t^{2m+2}}{s(s+2) \cdots (s+2m)}
\int_0^\infty e^{- \pi t^2 x^2} x^{s+2m+1}dx & = & \frac{\pi^{m+1}
t^{-s}}{\frac{s}{2} \left( \frac{s}{2} +1 \right) \cdots \left(
\frac{s}{2} +m\right)}
\int_0^\infty e^{- \pi y^2} y^{s+2m+2} \frac{dy}{y} \nonumber \\
& = &
\frac{\pi^{m+1} t^{-s} \pi^{- \frac{s+2m+2}{2}}
\Gamma \left( \frac{s+2m+2}{2} \right)}
{\frac{s}{2} \left( \frac{s}{2} +1 \right) \cdots
\left( \frac{s}{2} +m\right)}
 \nonumber \\
& = & \pi^{- \frac{s}{2}} \Gamma \left( \frac{s}{2} \right) t^{-s} \,,
\end{eqnarray}
where we repeatedly used $s \Gamma (s) = \Gamma (s+1).$ The
definition \eqn{P310} is valid for $0 \le t < \infty$ and defines
$h(t,s)$ as a meromorphic function for $\Re (s+2m+1) > 0$, and
\eqn{P311} shows that the definition of $h(t,s)$ is independent of
$m$. It follows that $h(t,s)$ is meromorphic for $s \in \CC$ with
possible poles at $s \in 2\ZZ_{\le 0}$.
For fixed $s \not\in 2\ZZ_{\le 0}$ it is
continuous in $t$, and satisfies \beql{P312} h(0,s) \equiv 0 ~. \eeq
Now \eqn{P35} and $\lim\limits_{t \to 0^+} h (t,s) =0$ follow from
\eqn{P39}--\eqn{P312}.

We obtain the second formula \eqn{P37}
from  \eqn{P35} by multiplying by $t$, replacing $s$ with $s + 1$
and setting
$$k(t,s) := - \frac{t e^{- \pi t^2}}{s+1} + t~ h(t, s+1) \,$$
as the remainder term in \eqn{P37}. This function is continuous in
$t$ for $0 \le t < \infty$ and $k(0,s) \equiv 0$, except possibly
for $s \in 1+2\ZZ_{\le -1}$.
~~~$\bsq$\\

\paragraph{\bf Proof of Theorem \ref{th31}.}
We define 
\begin{eqnarray}\label{P313}
\hat{L}_0^T (s,a,c) & := & \int_1^\infty \left(\sum_{n \in \ZZ \atop n \neq 0,-1}
 e^{2 \pi ian} e^{- \pi (n+c)^2 x^2} \right) x^{s-1} dx
\nonumber \\
&& + \int_1^\infty \left( \sum_{n \in \ZZ \atop n \neq 0,1}
e^{2 \pi i c(n-a)} e^{- \pi (n-a)^2 y^2} \right) y^{-s} dy
\end{eqnarray}
and
\begin{eqnarray}\label{P314}
\hat{L}_1^T (s,a,c) & := & \int_1^\infty \sum_{n \in \ZZ \atop n \neq 0,-1}
e^{2 \pi ian} (n+c) e^{- \pi (n+c)^2 x^2} x^s dx \nonumber \\
&& - i \int_1^\infty \left(
\sum_{n \in \ZZ \atop n \neq 0,1} e^{2 \pi ic (n-a)} (n-a)
e^{- \pi (n-a)^2 y^2} \right) y^{1-s} dy \,.
\end{eqnarray}
It is easy to see that $\hat{L}_0^T (s,a,c)$ and $\hat{L}_1^{T}(s,a,c)$ are
entire functions in $s$ and are jointly continuous in $(a,c)$ for
$(a,c) \in (-1,2) \times (-1,2)$.
For $0 < a < 1$ and $0 < c < 1$ we have, by definition,
\begin{eqnarray}\label{P315}
&&\hat{L}^{+}( s,a,c) = \hat{L}_0^T(s,a,c) +
\int_1^\infty e^{- \pi c^2 x^2} x^{s-1} dx +
e^{-2 \pi ia} \int_1^\infty e^{- \pi (c-1)^2 x^2}
x^{s-1} dx \nonumber \\
&&+ \int_1^\infty e^{-2 \pi iac}
e^{- \pi a^2 y^2} y^{-s} dy +
\int_1^\infty e^{2 \pi ic (1-a)}
e^{- \pi (1-a)^2 y^2} y^{-s} dy,
\end{eqnarray}
and
\begin{eqnarray}\label{P316}
&&\hat{L}^{-}(s,a,c) =  \hat{L}_1^T (s,a,c) + \int_1^\infty
ce^{- \pi c^2 x^2} x^s dx +
\int_1^\infty e^{-2 \pi ia} (c-1) e^{- \pi (c-1)^2 x^2} x^s dx \nonumber \\
&& - i \int_1^\infty e^{-2 \pi iac} (-a)e^{- \pi a^2 y^2} y^{1-s} dy
 - i \int_1^\infty
e^{2 \pi ic (1-a)} (1-a) e^{- \pi (1-a)^2 y^2} y^{1-s} dy \,.
\end{eqnarray}
We evaluate $\zeta_\ast (s,a,c)$ using the formula \beql{P317}
\zeta_\ast (s,a,c) = \frac{1}{2} \left( \frac{
\hat{L}^{+}(s,a,c)}{\pi^{-\frac {s}{2}} \Gamma (s/2)} + \frac{
\hat{L}^{-}( s,a,c)}{\pi^{- \frac{1+s}{2}} \Gamma \left(
\frac{1+s}{2} \right)} \right) \,. \eeq When considering the limit
of $  \hat{L}^{\pm}(s,a,c)$ as $(a, c)$ approaches a point $(a',
c')$ on the boundary of $\Box$, only those integrals on the right
side of \eqn{P315} (resp. \eqn{P316}) in which the exponential
function has exponent approaching 0 can introduce discontinuities.
We call these the {\em bad terms}. All the remaining terms
converge to the corresponding terms in $ \hat{L}^{\pm}(s, a',c' )$
in \eqn{JP240} and are represented by the corresponding integrals.
At the point $(a',c')$ the bad term in $\zeta_{*}(s, a', c')$ is
 replaced by the
quantity \beql{P318} \delta_{\ZZ}(a') \frac{1}{s-1} -
\delta_{\ZZ}(c') \frac{e^{- 2 \pi i a' c'}} {s}. \eeq The problem
is that for certain $s$ the  bad terms do not always continuously
approach these values as $(a,c) \to (a', c')$.

\paragraph{\em Case I (c=1)}
 $(a', c') = (a', 1)$ with $0 < a' < 1$.

There is exactly one bad term in $\hat{L}^{+}( s,a,c)$,
namely $e^{-2 \pi ia} \int_1^\infty e^{- \pi (1-c)^2 x^2} x^{s-1}
dx$. For it Lemma \ref{le31} gives, for $s \in \CC \smallsetminus 2\ZZ_{\le 0}$,
\beql{P319}
 e^{-2 \pi ia} \int_1^\infty e^{- \pi (1-c)^2 x^2}
x^{s-1} dx = -e^{-2 \pi ia} \frac{e^{-\pi (1-c)^2}}{s} + \pi^{-
\frac{s}{2}} \Gamma \left( \frac{s}{2} \right) e^{-2 \pi ia}
(1-c)^{-s} + e^{-2 \pi ia} h(1-c, s),
 \eeq 
 with $h(1-c,s) \to 0$ as $c \to 0$. 
 As $(a,c) \to (a',1)$, the first term on the right side converges to 
 $\frac{-e^{-2 \pi ia'}}{s}$, and the last term converges to zero.

There is one bad term in $\hat{L}^{-}(s,a,c)$, and for it Lemma
\ref{le31} gives, for 
$s \in \CC \smallsetminus \{ -1+2\ZZ_{<0}\}$,

\beql{P320}
 e^{- 2 \pi ia} \int_1^\infty (c-1) e^{- \pi (c-1)^2 x^2}
x^s dx = - \pi^{- \frac{s+1}{2}} \Gamma \left( \frac{s+1}{2} \right)
e^{- 2 \pi ia} (1-c)^{-s} + e^{- 2 \pi ia} k(1-c,s) \,, 
\eeq
 where
$k(1-c,s) \to 0$
as  $c \to  1$. When
these formulas are substituted in \eqn{P317}, the middle term on the
right side of \eqn{P319} cancels with the gamma-factor term on the
right side of \eqn{P320}. This yields \eqn{P31}, when 
$s \in \CC \smallsetminus \ZZ_{\le 0}$.

\paragraph{\em Case II. (c=0)}
$(a', c' ) = (a', 0)$ with $0 < a' < 1$.

Here $\hat{L}^{+}( s,a,c)$ has a single bad term, for which
Lemma \ref{le31} gives, for $s \in \CC \smallsetminus 2\ZZ_{\le0}$,
\beql{P321}
\int_1^\infty e^{- \pi c^2 x^2} x^{s-1} dx = -
\frac{e^{- \pi c^2}}{s} + \pi^{- \frac{s}{2}} \Gamma
 \left( \frac{s}{2} \right) c^{-s} + h(c,s) \,.
\eeq
As $(a,c) \to (a', 0)$
the first term on the right converges to the constant $- \frac{1}{s}$, and the
last term converges to zero. The middle term, when substituted in
\eqn{P317}, contributes $\frac{1}{2} c^{-s}$ to the right side of \eqn{P317}.

Now $\hat{L}^{-}( s,a,c)$ has a single bad term, and a
similar analysis shows, for $s \in \CC \smallsetminus\{-1+2\ZZ_{\le0}\}$,
that it also contributes $\frac{1}{2} c^{-s}$ to the right side of \eqn{P317}. Put
together, we obtain as $(a,c) \to (a', 0^+),$
\begin{eqnarray*}
\zeta_\ast (s,a,c) - c^{-s} & = & \frac{1}{2}
\left( \frac{\hat{L}^{+}(s,a,c)}{\pi^{- \frac{s}{2}}
\Gamma \left( \frac{s}{2} \right)} -  c^{-s} \right)
+ \frac{1}{2}
\left( \frac{\hat{L}^{-}(s,a,c)}{\pi^{- \frac{s+1}{2}}
\Gamma \left( \frac{s+1}{2}\right)} -  c^{-s} \right) \nonumber \\
& \to &\zeta_\ast (s, a',0),
\end{eqnarray*}
which is \eqn{P32}, when $s \in \CC \smallsetminus \ZZ_{\le 0}$.

\paragraph{\em Case III (a=1)} 
$(a', c') = (1,c')$, with $0 < c' < 1$.

There is one bad term in $\hat{L}^{+}(s,a,c)$, which using
Lemma \ref{le31} is, for $1-s \in \CC \smallsetminus 2\ZZ_{\le 0}$,
\begin{eqnarray*}
&&
 \int_1^\infty e^{2 \pi ic (1-a)} e^{- \pi (1-a)^2 y^2} y^{-s}dy =
- e^{2 \pi ic (1-a)}
\frac{e^{-\pi (1-a)^2}}{1-s} \\
&&+ \pi^{- \frac{1-s}{2}} \Gamma \left( \frac{1-s}{2} \right)
e^{2 \pi i c(1- a)} (1-a)^{s-1} + e^{2 \pi i c(1 - a)} h(1-a, 1-s) \,.
\end{eqnarray*}
As $(a, c) \to (1, c' )$ the first term on the right side converges
to the constant $\frac{1}{s-1}$,  and the last term vanishes. 
There is one bad term in
$\hat{L}^{-}(s,a,c)$, which using Lemma \ref{le31} is, for
$1-s \in \CC \smallsetminus \{-1+2\ZZ_{\le 0}\}$,
\begin{eqnarray*}
 - i \int_1^\infty e^{2 \pi ic (1-a)} (1-a) e^{- \pi (1-a)^2 y^2}
y^{1-s} dy = &- &i \pi^{- \frac{2-s}{2}} \Gamma \left(
\frac{2-s}{2} \right) e^{2 \pi ic (1-a)}
(1-a)^{s-1} \\
&- &ie^{2 \pi ic (1-a)} k(1-a, 1-s) \,,
\end{eqnarray*}
in which the second term on the right side vanishes as 
$(a,c) \to(1, c')$. Substituting these into \eqn{P317}, and subtracting off
the remaining two terms of the right sides yields, as $(a,c) \to
(1,c')$, 
\beql{P322} 
\zeta_\ast (s,a,c) - \ggG(s) e^{2 \pi ic(1-a)} (1-a)^{s-1}~~ \to~~ \zeta_\ast (s,1,c'),
 \eeq
 in which
 \beql{P323} 
   \ggG(s) :=
\frac{1}{2} \left( \frac{\pi^{- \frac{1-s}{2}} \Gamma \left(
\frac{1-s}{2} \right)}{\pi^{- \frac{s}{2}} \Gamma \left(
\frac{s}{2} \right)} - i \frac{\pi^{- \frac{2-s}{2}} \Gamma \left(
\frac{2-s}{2} \right)} {\pi^{- \frac{s+1}{2}} \Gamma
\left(\frac{s+1}{2} \right)} \right) \,.
 \eeq 
 Using the identity
$\frac{\Gamma ( \frac{2-s}{2} )}{\Gamma \left( \frac{s+1}{2}
\right)} = \frac{\Gamma \left( \frac{1-s}{2}
\right)}{\Gamma (\frac{s}{2})} \cot \frac{\pi s}{2} $,
 the formulas \eqn{P322} and \eqn{P323} yield \eqn{P33}, for
 $s \in \CC \smallsetminus \ZZ_{\ge1}$.

\paragraph{\em Case IV (a=0)}
$(a', c') = (0, c')$, with $0 < c' < 1$.

This case is  parallel to case III.~~~$\bsq$ \\

It is clear from the proof above that the case when $(a,c)$ tends to
a corner of the unit square can be handled in a similar fashion,
except that each of $\hat{L}^{+}(s,a,c)$ and
$\hat{L}^{-}(s,a,c)$ now  have two bad terms which must be
taken into account.  The result for this case is
as follows.

\begin{theorem}\label{th32}
{\em (Lerch Corner Limits)}
For fixed $s \in \CC \smallsetminus \ZZ$ and $(a,c) \in \Box^\circ$ the
following hold.

\begin{itemize}
\item[(i)]
As $(a,c) \to (1^- , 1^- )$,
\begin{eqnarray}~\label{P324}
\zeta_\ast (s,a,c) & - & \frac{1}{2}
\frac{\pi^{- \frac{1-s}{2}}
\Gamma ( \frac{1-s}{2} )}{\pi^{- \frac{s}{2}} \Gamma ( \frac{s}{2} )}
\left( 1- i \cot \frac{\pi s}{2} \right)
e^{2 \pi i(1 - a)c} (1-a)^{s-1} \to
\zeta_\ast (s,1,1)  \,. \nonumber
\end{eqnarray}

\item[(ii)]
As $(a,c) \to (0^+ , 1^-)$,
\begin{eqnarray}\label{P325}
\zeta_\ast (s,a,c) & - & \frac{1}{2}
\frac{\pi^{- \frac{1-s}{2}}
\Gamma (\frac{1-s}{2} )}{\pi^{- \frac{s}{2}} \Gamma ( \frac{s}{2} )}
\left( 1+ i \cot \frac{\pi s}{2} \right) e^{- 2 \pi iac} a^{s-1} \to
\zeta_\ast (s,0,1) \,. \nonumber
\end{eqnarray}

\item[(iii)]
As $(a,c) \to (1^- , 0^+)$,
\begin{eqnarray}\label{P326}
 \zeta_\ast (s,a,c)& - & \left(c^{-s} +
\frac{1}{2} \frac{\pi^{- \frac{1-s}{2}} \Gamma ( \frac{1-s}{2})}
{\pi^{- \frac{s}{2}} \Gamma ( \frac{s}{2} )} \left( 1- i \cot
\frac{\pi s}{2} \right) e^{2 \pi i(1-a)c} (1-a)^{s-1} \right) ~~
\to~~ \zeta_\ast (s,1,0) \,. \nonumber
\end{eqnarray}

\item[(iv)]
As $(a,c) \to (0^+ , 0^+ )$,
\begin{eqnarray}\label{P327}
\zeta_\ast (s,a,c) &- &\left( c^{-s} + \frac{1}{2}
\frac{\pi^{\frac{1-s}{2}}
\Gamma ( \frac{1-s}{2} )}{\pi^{- \frac{s}{2}} \Gamma ( \frac{s}{2} )}
\left( 1+ i \cot \frac{\pi s}{2} \right) e^{-2 \pi iac} a^{s-1} \right)
~~ \to~~ \zeta_\ast (s,0,0) \,. \nonumber
\end{eqnarray}

\end{itemize}
\end{theorem}

\pf
The details  are similar to those of Theorem~\ref{th31}; we
omit them. 
Note that the limits on the right in cases (i) - (iv)
are all $\zeta(s)$, since
$$
\zeta_\ast (s,0,0) =\zeta_\ast (s,0,1)= \zeta_\ast (s,1,0)
=\zeta_\ast (s,1,1)= \zeta(s).
$$
The limits above
are also valid when  $s$ is a negative
odd integer. ~~~$\bsq$ \\

We now deduce Theorem \ref{th14} using Theorem~\ref{th31} and \ref{th32}.\\

\paragraph{\bf Proof of Theorem \ref{th14}.}
We apply Theorems \ref{th31} and \ref{th32} for the case
$s \in \CC \smallsetminus  \ZZ$,
and treat remaining case $s \in \ZZ$ separately.

Suppose that $s \in \CC \smallsetminus \ZZ$.
 Theorem \ref{th31}(i) gives a continuous
extension to the edge $c=1$, $0 < a' < 1$ for all $s\in \CC$.
Theorem \ref{th31}(ii) gives a continuous extension to the edge
$c=0$, $0 < a' < 1$ exactly when $c^{-s}$ continuously extends,
which is when $\Re (s) < 0.$ Theorem \ref{th31}(iii) gives a
continuous extension to the edge $a=1$, $0 < c' < 1$ exactly when
$(1-a)^{s-1}$ continuously extends, which is when $\Re (s) > 1$.
Theorem \ref{th31}(iv) gives a continuous extension to the edge
$a=0$, $0 < c' < 1$ exactly when $a^{s-1}$ continuously extends,
which is when $\Re (s) > 1$. Theorem \ref{th32} covers the four
corners similarly giving a continuous extension to the two corners
$(a', c') = (0,1)$ and $(1,1)$ when $\Re (s) > 1$ and in no other
cases.

Now suppose that $s= m \in \ZZ$.
If $m \ge 2$ then we are in the region of absolute convergence
of the Dirichlet series
\beql{P334}
\zeta (s,a,c) = c^{-s} + \sum_{n=1}^\infty e^{2 \pi ina} (n+c)^{-s} \,,
\eeq
so we get a continuous extension of the sum in the right side of
\eqn{P334} to $\Box$.
The term $c^{-s}$ produces a discontinuity on the side $c=0$,
including the two corners.
The case $m \le -1$ is handled using Lerch's transformation formula
\eqn{sp34}, setting $1-s = m$.
Then $s= |m|+1$ in \eqn{sp34} so the right side has
absolutely convergent Dirichlet series, which are continuous
except for two terms producing
discontinuities along the lines $a=0$, $a=1$ and the four corners.

Next suppose that $s=1$.
When $a'=0$ or 1 the function $\zeta_\ast (s,a',c)$ has a pole at $s=1$,
so $\zeta_\ast (s, a,c)$ diverges as $(a,c) \to (0,c')$ or $(1, c')$.
For $s=1$ and $(a,c) \in \Box^\circ$ the sum in \eqn{P334}
conditionally converges, uniformly in $\ep < a < 1- \ep$, while
the term $c^{-s}$ diverges as $c \to 0^+$.
Thus there is no continuous extension to any point on the line $c' =0$.
There is however a continuous extension to the line $c' =1$, $0 < a' < 1$,
by the same reasoning.

The remaining case is $s=0$. The Lerch transformation formula
\eqn{sp34} for $s=1$ expresses $\zeta_\ast (0,a,c)$ on the lines
$a=0$ and $a=1$ as a linear combination of $\zeta_\ast (1,1-c,a)$
and $\zeta_\ast (1,c,1-a)$, with both coefficients nonzero. One of
these has a continuous extension to the boundary and the other
doesn't, so their sum never does. We conclude from
Theorem~\ref{th31}(i),(ii) and the holomorphicity of $\zeta_*(s, a,
c)$ at $s=0$ that $\zeta_*(0, a, c)$ has a continuous extension to
the line $c = 1$, but not to the line $c = 0$
and not at the four corners of the square. ~~~$\bsq$ \\

%
%
%
%
%

\section{Renormalized Lerch functions}
\setcounter{equation}{0}

Theorem \ref{th31} and Theorem \ref{th32} provide a way to restore
continuity on the boundary of the unit square $\Box$ for the functions
$L^{\pm}(s,a,c)$ by ``renormalizing'' the functions with
correction terms,
 for $s \in \CC \smallsetminus \ZZ$, as follows. This renormalization
 is most easily expressed using  the following functions. \\

\begin{defi}~\label{de60}
{\em 
The {\em Tate gamma functions} $\gamma^{\pm}(s)$, 
sometimes called {\em Gelfand-Graev gamma functions}, are given by}
$$
\ggG^{+}(s)   := \frac{ \pi^{-\frac{s}{2}} \Gamma( \frac{s}{2})}
{ \pi^{-\frac{1-s}{2}} \Gamma(\frac{1-s}{2})},  \,\,\,\,\,\,\,\,\,
\ggG^{-}(s)   :=  \frac{ \pi^{-\frac{s+1}{2}} \Gamma( \frac{s+1}{2})}
{ \pi^{-\frac{2-s}{2}} \Gamma(\frac{2-s}{2})}, 
$$
\end{defi}

These functions were introduced in
Tate's thesis \cite[p. 317]{Ta67}  as  normalizing factors in local  functional equations
at the real place,
see also Burnol \cite[p. 822]{Bu04}.
They satisy  the identities
\beql{N114}
\ggG^{\pm}(s) \ggG^{\pm}(1-s) = 1 ~~~\mbox{for}~~ s \in \CC.
\eeq
The functions $\ggG^{\pm}(1-s)= \ggG^{\pm}(s)^{-1}$ have  a "scattering matrix" interpretation, 
given in Burnol \cite[Sec. 5]{Bu03}.
These functions  are present in various formulas in 
\S5, where \eqn{P323}  can be rewritten 
\beql{N114b}
\gamma(s) = \frac{1}{2}\left( \gamma^{+}(1-s) -i \gamma^{-}(1-s) \right),
\eeq
and  the functions  $\gamma^{\pm}(1-s)$ appear in formulas  of Theorem~\ref{th32}.


\begin{defi}~\label{de51} 
{\em Let $s \in \CC \smallsetminus \ZZ$,
and $k = 0$ or $1$. Write  $(-1)^k= \pm $ and define for  $(a,c) \in \Box^{\circ}$ the
{\em renormalized Lerch functions}
 \beql{P328} 
 {L}^{R,\pm}( s,a,c) := L^{\pm}( s,a,c) - S^{\pm}( s,a,c),
\eeq
obtained from $L^{\pm}(s, a, c)$ by removing the "correction term"
\begin{eqnarray}\label{P328b}
S^{\pm}( s,a,c) & := &  
\{ c^{-s} +(-1)^k e^{-2 \pi ia} (1-c)^{-s} \} \nonumber \\
& + & i^k \ggG^{\pm}(1-s) \{ e^{-2 \pi iac} a^{s-1} +(-1)^k e^{2\pi i(1- a)c} (1-a)^{s-1} \}   ,
\end{eqnarray}
in which $\ggG^{\pm}(s)$ denotes the Tate gamma function
given in  definition \ref{de60}.
 }
\end{defi}

The function $S^{\pm}( s,a,c)$ embodies the total
contribution of the ``bad terms''. This definition cannot be further
extended to $s \in \ZZ$, because for each integer $s$ at least one
of the gamma factors in  $S^{\pm}( s,a,c)$ has a pole. However this
definition does further extend in the $(a,c)$-variables to
$(a,c) \in (\RR \smallsetminus \ZZ) \times (\RR \smallsetminus \ZZ)$.

\begin{theorem}\label{th31A}
{\em (Renormalized Lerch Zeta Functions)}

(1) For fixed $s \in \CC \smallsetminus \ZZ$
the ``renormalized'' functions ${L}^{R,\pm}(s,a,c)$
continuously extend to the closed unit square $\Box$ in the $(a,c)$ variables. 

(2) For fixed $(a, c) \in \Box$ the ``renormalized " functions ${L}^{R,\pm}(s,a,c)$
are meromorphic functions in $s$ having at most simple poles, with all poles 
 in the set $\ZZ$. 

(3) For all $s \in \CC \smallsetminus \ZZ$
and $(a,c) \in \Box$ the completed  ``renormalized" functions
$$
\hat{L}^{R, \pm}( s,a,c) := 
\pi^{-(\frac{s+k}{2})}\Gamma\left(\frac{s+k}{2}\right)
 \hat{L}^{R, \pm}(s, a, c)
$$
satisfy the functional equations
\beql{P330}
\hat{{L}}^{R,\pm}(s,a,c) = w_{\pm} e^{- 2 \pi iac}
 \hat{{L}}^{R,\pm}( 1-s, 1-c, a), 
\eeq
with $w_{+}=1, w_{-}=i$.
\end{theorem}

\paragraph{\bf Remark.} This effect of this ``renormalization'' is to remove the
singularities on the boundary of the unit square. However  the ``renormalized"
functions ${L}^{R, \pm }(s, a, c),$  when extended to
$(a,c) \in (\RR \smallsetminus \ZZ) \times (\RR \smallsetminus \ZZ)$, still have
singularities at all integer values of $a$ or $c$ other than $0$ or
$1$.\\

\pf
For each fixed $(a, c) \in \Box^{\circ}$
the functions $S^{\pm}( s, a, c)$ are
meromorphic functions of $s$,
with simple poles confined to $s \in \ZZ$. Theorem~\ref{th12}
states that
$\hat{L}^{\pm}( s,a,c)$ also has this  property, with
pole set confined to $s \in \{0, 1\}$. Thus
${L}^{R,\pm}( s,a,c)$ inherits this property
for $ (a, c) \in \Box^{\circ}$.

The functional equation \eqn{P330} for $(a, c) \in \Box^{\circ}$ is
inherited from the fact that $\hat{L}^{\pm}(s,a,c)$
satisfies this functional equation (Theorem~\ref{th12}), as does the
completed ``renormalizing'' function
 $\hat{S}^{\pm}( s, a, c):= \pi^{-(\frac{s+k}{2})}\Gamma\left(\frac{s+k}{2}\right)S^{\pm}( s, a, c)$.
 Indeed one easily verifies  the relation
\beql{P331aa}
 \hat{S}^{\pm}(s, a, c) = i^k e^{-2\pi i a c}
\hat{S}^{\pm}( 1-s, 1-c, a).
\eeq

We now establish that  a continuous extension to the boundary of $\Box$
exists for the case 
${L}^{R,+}(s,a,c)$; we omit details for the case 
${L}^{R,-}(s,a,c)$, which is similar.
In terms of the expression in Theorems \ref{th31} and \ref{th32}, we
can write
\beql{P331}
{L}^{R,+}( s,a,c) = A(s,a,c) + B(s,a,c) \,,
\eeq
in which
\begin{eqnarray}\label{P332}
&A(s,a,c)\quad := \quad  \zeta (s, a, c) - c^{-s} \qquad \qquad \qquad
 \qquad \qquad \qquad \qquad \qquad \qquad \qquad \qquad \nonumber \\
&-  \frac{1}{2}\ggG^{+} (1-s) \{ ( 1 + i \cot \frac{\pi s}{2} )
e^{-2 \pi iac} a^{s-1} + ( 1- i \cot \frac{\pi s}{2}) e^{2 \pi
i(1-a)c} (1-a)^{s-1} \},
\end{eqnarray}
\begin{eqnarray}\label{P333}
&B(s,a,c)\quad := \quad  e^{-2 \pi i a}\zeta (s, 1-a,1-c) - e^{- 2
\pi ia} (1-c)^{-s} \qquad \qquad \qquad
 \qquad \qquad \quad  \nonumber \\
&- \frac{1}{2} \ggG^{+} (1-s) \{ ( 1- i \cot \frac{\pi s}{2})
e^{-2 \pi iac} a^{s-1} + ( 1+ i \cot \frac{\pi s}{2} ) e^{2 \pi
i(1-a)c} (1-a)^{s-1} \}\,
\end{eqnarray}
where
$\ggG^{+} (s) = \frac{\pi^{- \frac{s}{2}} \Gamma \left(\frac{s}{2} \right)}
{\pi^{- \frac{1-s}{2}} \Gamma \left( \frac{s1-}{2}\right)}$.
Since $s \not\in \ZZ$ none of the $\Gamma$-factors or
$\cot \frac{\pi s}{2}$ have poles. In addition 
$\zeta (s,a, c) =\zeta_\ast (s,a,c)$ for $(a,c) \in \Box^\circ$,
so that Theorem \ref{th31} applies. We claim that $A(s,a,c)$ continuously extends to
the boundary for $s \in \CC \smallsetminus \ZZ$.
To see this, when approaching the
line $c=1$, $0 < a' < 1$,
the function $\zeta_\ast(1,a,c)$ continuously extends to it 
by Theorem \ref{th31}(i), and the
other three terms on the right side of \eqn{P332} also continuously
extend there. When
approaching the line
$c=0$, $0 < a' < 1$, the term $\zeta (s,a,c) -c^{-s}$ continuously
extends by Theorem \ref{th31}(ii), while the other two terms on the
right side of \eqn{P332} are continuous there. When
 approaching the line $a=1$, $0 < c' < 1$, the term
$\zeta_\ast (s,a,c) - \frac{1}{2} \ggG^+ (1-s) (1- i \cot\frac{\pi s}{2} ) e^{2 \pi i(1 -a)c} (1 -a)^{s-1}$
 has a continuous extension by Theorem \ref{th31}(iii) and the other two
terms are continuous there. Approaching the line 
$a=0$, $0 < c' <1$, the term 
$\zeta_\ast (s,a,c) - \frac{1}{2} \ggG^+ (1-s) (1+ i\cot \frac{\pi s}{2} ) e^{-2 \pi iac} a^{s-1}$ 
has a continuous
extension by Theorem \ref{th31}(iv) and the other two terms are
continuous there.
To get continuity at the four corners of $\Box$, Theorem \ref{th32} is
invoked in a similar fashion. This proves the claim.

In a similar fashion one proves that $B(s,a,c)$ continuously extends
to the boundary of $\Box$ when $s \in \CC \smallsetminus \ZZ$.
Now
\eqn{P331} gives a continuous extension of
${L}^{R,+}(s,a,c)$ to $\Box$.

The continuous extension to $\Box$ implies that the functional
equation \eqn{P330}  continues to hold on the boundary
 $(a, c) \in \partial \Box = \Box \smallsetminus \Box^{\circ}$.
 It remains to establish that  for
$(a, c) \in \partial \Box$  the extended functions
${L}^{R,\pm}( s, a, c)$
 are meromorphic functions of $s \in \CC$, having at most simple poles
for $s \in \ZZ$. The crucial fact needed to do this  is that the difference
functions obtained from subtracting off the ``bad terms'' are meromorphic in
$s$, which is the content of Lemma~\ref{le31};
we omit  details. ~~~$\bsq$ \\

We now apply Theorem~\ref{th31A} to deduce Theorem~\ref{th15},
which answers the question:  For which 
values 
$s \in \CC$ do the functions
$L^{\pm}( s, a, c)$ belong to $L^p(\Box, da~ dc)$, for fixed $1\le p \le 2$?\\

\paragraph{\bf Proof of Theorem ~\ref{th15}.}
We treat first the case of  $s \in \CC \smallsetminus \ZZ$.
Theorem~\ref{th31A} asserts that
the functions ${L}^{R,\pm}( s, a, c)$ are
 continuous on the unit
square (except for integral $s$), so the property of belonging to $L^p(\Box, da dc)$ 
for $L^{\pm}( s, a, c)$ is completely determined by the
behavior of  the ``correction'' functions $S^{\pm}( s, a,c)$, and
similarly for linear combinations 
$F_s(a, c) =c_1L^{+}( s, a, c)+ c_2 L^{-}(s, a, c) \in \sE_s$.
It then suffices to prove the result for the completed
functions $\hat{L}^{\pm}(s,a,c)$, because the gamma factors
take finite nonzero values there. 

Suppose first that  $ 1 \le p < 2.$
 In any linear
combination $F_s(a, c)$, one of the pair of correction terms
$c^{-s}$ and  $ e^{-2\pi i a}(1- c)^{-s}$ has a nonzero coefficient,
and one of the pair of  correction terms $e^{-2\pi i a c}a^{s-1},$ and $e^{2\pi i(1-a)c}(1-a)^{s-1}$
 has a nonzero coefficient. Furthermore the singularities of each of these terms
 cannot be cancelled by any linear combination of the other three terms.
 All functions in
both pairs belong to $L^p(\Box, da~dc)$ for $1 - \frac{1}{p}< \Re(s) < \frac{1}{p}$; thus
all linear combinations  of $L^{\pm}( s, a, c)$ belong to  $L^p(\Box, da~dc)$ in
this region.  Each of the first pair $c^{-s}$ and  $ e^{-2\pi i a}(1- c)^{-s}$ does not belong to 
$L^p(\Box,da~dc)$ 
for $\Re(s) \ge \frac{1}{p}$, and each of the second pair
$e^{-2\pi i a c}a^{s-1},$ and $e^{2\pi i(1-a)c}(1-a)^{s-1}$
does not belong
to $L^{p}(\Box, da dc)$  for
$\Re(s) \le 1- \frac{1}{p}$, hence all nonzero linear combinations $F_s(a, c)$ do not belong
to $L^p(\Box, da~dc)$ for  $\Re(s) \le 1 - \frac{1}{p}$ and for $\Re(s) \ge \frac{1}{p},$
for non-integral $s$.

To treat the case $p\ge 2$, it suffices to show
the result for $p=2$, since $L^p(\Box, da dc) \subset L^2(\Box, da dc)$ if $p \ge 2$.
We note that  each of $c^{-s}$ and  $ e^{-2\pi i a}(1- c)^{-s}$ 
does not belong to
$L^2(\Box, da~ dc)$ for $\Re(s) \ge \frac{1}{2}$, and each of 
 $e^{-2\pi i a c}a^{s-1},$ and $e^{2\pi i(1-a)c}(1-a)^{s-1}$
does not belong to $L^2(\Box, da~ dc)$
 for $\Re(s) \le \frac{1}{2}$. As above, we conclude  that no nonzero
  $ F_{s}(a, c)$ belongs to  $L^2(\Box, da dc)$, for all $s \in \CC \smallsetminus \ZZ$.
Thus assertions (1) and (2) are verified in this case.

It remains to consider the case where $s \in \ZZ$. We now observe that the
Dirichlet series for $L^{\pm}( s, a, c) - |c|^{-s} -
(-1)^ke^{-2\pi i a} |1-c|^{-s}$ converges conditionally for
$\Re(s)>0$ and defines a continuous function on 
$\Box^{\circ} \cup\{(a,0): 0 < a <1\} \cup  \{(a,1): 0 < a <1\}$.
 As the edges $c=0$ and $c=1$ are approached, the
growth of the functions $ |c|^{-s}$ and $|1-c|^{-s}$ puts the
functions $L^{\pm}( s, a, c)$ not in $L^p( [\delta, 1-\delta]\times [0,1], da~dc)$
 for any $\delta>0$ and $\Re(s) \ge 1$.
This includes $s=1, 2,3,...,$ so $L^{\pm}( s, a, c) \not\in L^p(\Box, da~dc)$ for
$s=1,2,3...$ The same argument applies to exclude membership in
$L^2(\Box, da~dc)$ for $s=1,2,3...$. All nonzero linear combinations
of these functions are excluded similarly, since it is not possible
to cancel both terms $ |c|^{-s}$ and $|1-c|^{-s}$ in any linear
combination. Finally, a similar argument for $s=0, -1, -2, ...$
 excludes  all nonzero linear combinations of
$e^{-2\pi i ac}L^{\pm}( 1-s, 1-c, a)$ 
from membership in $L^p(\Box, da~dc)$ for $1 \le p \le 2$,
 since these functions have a Dirichlet
series expansion absolutely convergent in $\Re(s) < 0$
and conditionally convergent in $\Re(s) < 1$, inside $\Box^{\circ}$.
$~~~\bsq$

%
%
%
%
%

\section{Weyl algebra action}
\setcounter{equation}{0}

The Tate viewpoint on the functional equation uses test functions
drawn from the Schwartz space $\sS= \sS(\RR)$.
The Schwartz space is closed under
the action of the {\em Weyl algebra} 
$\bA_1:= \CC [ x, \frac {\partial}{\partial x}],$
and carries a (smooth) representation of the Weyl algebra  acting as linear operators on $\sS$.
The Mellin transform intertwines the Weyl algebra action with
an action on functions in the $s$-variable given by difference operators. We use this to define for
integer $n  \ge 0$ an
infinite family $\hzt_n(s, a, c)$ of generalized Lerch zeta functions
with functional equation, associated to the oscillator representation,
whose first two members for $n=0,1$ are $\hat{L}^{+}(s, a,c)$
and $\hat{L}^{-}(s, a, c)$, respectively.

The {\em Weyl algebra} is the universal enveloping algebra $U( \bh_{\RR})$
of the real Heisenberg Lie algebra $\bh_{\RR}.$ Here the {\em real
Heisenberg Lie algebra} 
is the  three-dimensional real Lie algebra
\beql{Z501}
 \bh_{\RR} := \RR[I, x, \frac{\partial}{\partial x}] 
= \RR [I, D_+, D_-] \ , 
\eeq 
in which
 \beql{405}
  D_+ := \sqrt{2 \pi}
\left( x - \frac{1}{2 \pi} \frac{\partial}{\partial x} \right)
\quad\mbox{and}\quad D_- := \sqrt{2 \pi} \left( x + \frac{1}{2 \pi}
\frac{\partial}{\partial x} \right) ~. 
\eeq 
Its commutation
relations on the generators $I$, $D_+$, $D_-$ are 
\beql{406}
 [D_+,D_- ] = -2I \quad\mbox{and} \quad [D_+ , I ] =[  D_- , I ] =0 ~.
\eeq

The Mellin transforms $\sM_{k}$  acting on even and odd
Schwartz functions serve as intertwining operators to carry the
Weyl algebra action  on functions on $\RR^{\ast}$
over to function spaces in the $s$-variable.
\begin{lemma}\label{nle41}
The Weyl algebra $\bA_1$ acts on
the two-sided Mellin transforms $\sM_k$ of
$f \in \sS$,  with $k=0, 1~ (\bmod ~2)$, as
difference operators in the $s$-variable, as follows:
\begin{eqnarray}
\sM_k (\frac{\partial f}{\partial x})(s)  & = &
-(s - 1) \sM_{k+1} (f)( s-1) ~,\label{402}\\
&& \nonumber \\
\sM_k (xf(x)) (s) & = & \sM_{k+1} (f)(s+1) ~.\label{403}
\end{eqnarray}
As a consequence,
\begin{eqnarray}
\sM_k (\frac{\partial^2 f}{\partial x^2})(s) &  =&
(s - 1)(s-2) \sM_{k} (f)( s-2) ~, \label{801a}\\
&& \nonumber \\
\sM_k( x \frac{\partial f}{\partial x})(s) & = & -s \sM_{k} (f)(s)~ , \label{801b}\\
&& \nonumber\\
\sM_k (x^2 f(x)) (s)  & =  & \sM_{k} (f)(s+2) ~. \label{801c}
\end{eqnarray}

\end{lemma}

\pf
Integration by parts gives
\begin{eqnarray}\label{404}
\sM_k (\frac{\partial f}{\partial x})(s)  & = &
\int_{-\infty}^\infty
\frac{\partial f}{\partial x} (sgn\, x)^k |x|^s \frac{dx}{|x|} \nonumber \\
& = & -(s - 1) \int_{-\infty}^\infty
f(x) ( sgn \, x )^{k+1} |x|^{s-1} \frac{dx}{|x|} \nonumber \\
& = & -(s - 1) \sM_{k+1} (f)(s-1) ~.
\end{eqnarray}
The  calculation  \eqn{403} is immediate.
The remaining three identities follow by substitution.~~~$\bsq$ \\

We consider the Weyl algebra action on the eigenfunctions of the
one-dimensional harmonic oscillator Hamiltonian
$$\
\he{H} :=  \frac{1}{2}\left(D_{+} D_{-}+ D_{-}D_{+}\right)
=  -\frac {1}{2 \pi} \frac {\partial^2}{\partial x^2} + 2 \pi x^2.
$$
We work in the infinite-dimensional vector space
$$\sS_0 = \CC[x] e^{-\pi x^2}  \subset \sS, $$
which is closed under the action of the Weyl algebra.  Here
\beql{407a}
\sS_0 = \bigoplus_{n \ge 0} \CC[\phi_n],
\eeq
where the $\phi_n$ are a basis of
eigenfunctions of $\he{H}$,  given by $\phi_0(x)= e^{- \pi x^2}$ and 
\beql{407}
\phi_n (x) := D_+^n (\phi_0 (x)) = D_+^n (e^{-\pi x^2}),
\eeq
The functions $\phi_n(x)$ are  explicitly given by
\beql{408}
\phi_n (x) = 2^{- n/2} H_n ( \sqrt{2 \pi} x) e^{- \pi x^2},
\eeq
in which
$H_n (x)$ are the Hermite polynomials
\beql{409}
H_n (x) = (-1)^n e^{x^2} \frac{d^n}{dx^n} (e^{-x^2} ) ~.
\eeq
The functions $\phi_n(x)$ satisfy the (time-independent) Schr\"{o}dinger equation
\beql{410}
\he{H} \phi_n = (- \frac{1}{2 \pi}\frac{d^2}{dx^2}   +2 \pi x^2) \phi_n
 = (2n+1) \phi_n ~.
\eeq
The operators $D_{+}$ and $D_{-}$ are called {\em raising} and {\em lowering
operators}, due to the relations
$D_{+} \phi_n= \phi_{n+1}, D_{-}\phi_n = 2n\phi_{n-1}$.
The Fourier transform \eqn{203} acts compatibly with the raising operator
\beql{Z513}
\sF(D_+ f) = - i D_+ (\sF f) , \quad\mbox{all}\quad
f \in \sS ~,
\eeq
cf. Bump and Ng \cite[p. 198]{BN86}. (Bump and Ng use the convention
$\sF f (y) = \int_{-\infty}^\infty f(x) e^{2 \pi i xy} dx$,
so that \eqn{Z513} differs by  a sign change from their equation (1.4).)
Since $\phi_0 = \sF \phi_0$, we obtain from \eqn{407} by induction
on $n$ that
\beql{Z514}
\sF \phi_n (x) = (-i)^n \phi_n (x) ~, \quad n \ge 0 ~.
\eeq

Bump, Choi, Kurlberg, and Vaaler \cite{BCKV99} prove the following
results.
\begin{prop}\label{nle42}
For $n \ge 0$ the Mellin transforms of $\phi_n(x)$ satisfy
\begin{eqnarray}
\sM_0 (\phi_{2n}(x))(s) & = &
\pi^{- \frac{s}{2}} \Gamma \left( \frac{s}{2} \right) p_n (s), \label{411}\\
~~\sM_1 (\phi_{2n+1}(x))(s) & = &\sqrt{2\pi} \,
\pi^{- \frac{s+1}{2}} \Gamma \left( \frac{s+1}{2} \right) q_n (s),
\label{412}
\end{eqnarray}
in which $p_n (s)$ and $q_n (s)$ are real polynomials of degree $n$.
These polynomials satisfy the functional equations
\begin{eqnarray}\label{413}
p_n (1-s) & = & (-1)^n p_n (s) \,, \quad n \ge 0,
\nonumber \\
q_n (1-s) & = & (-1)^n q_n (s) \,, \quad n \ge 0.
\end{eqnarray}
The
zeros of $p_n (s)$, $q_n (s)$ all lie on the critical line $\Re (s) = 1/2$.
\end{prop}

\pf
This is proved in Bump et~al. \cite[Theorem 1]{BCKV99}.
The last assertion of this lemma is the ``local Riemann hypothesis''
studied in Bump and Ng~\cite{BN86},
Bump et al~\cite{BCKV99}, Kurlberg~\cite{K99}, and Olofsson\cite{Ol07}.~~~$\bsq$ \\

Lemma  ~\ref{nle41} yields   recursion relations for the
polynomials $p_n(s)$ and $q_n(s)$, namely
 \beql{413a}
\sM_k(\phi_{n+1})(s) =  \sM_k(D_+ \phi_n)(s)
 = 
\sqrt{2\pi} \sM_{k+1}(\phi_n)(s + 1) + \frac {1}{\sqrt{2\pi}} (s -
1) \sM_{k+1}(\phi_n)(s-1), 
\eeq 
where subscripts in $\sM_k$ are
given (mod ~$2$). From this one readily finds that the polynomials
$p_n$ and $q_n$ satisfy $p_0 (s) = q_0 (s) =1$ and
\begin{eqnarray}\label{414}
p_{n+1} (s) & = & s q_n (s+1) + (s - 1)q_n (s-1), \\
\label{415}
q_{n+1} (s) & = & p_{n+1} (s+1) + p_{n+1} (s-1) \,.
\end{eqnarray}
These in turn imply that  $p_n(x), q_n(x) \in \ZZ[x]$. The first
few polynomials $p_n(x)$ and $q_n(x)$ are given in Table 1.

\begin{table}[htb]
$$
\begin{array}{|c|r|r|} \hline
n & p_n(s)&  q_n(s) \\ \hline
0 & 1 & 1 \\
1&   2s-1   &   4s-2   \\
2 & 8s^2-8s+6 & 16s^2-16s+22 \\
3 &  32s^3-16s^2+92s-54&  64s^3-32s^2+376s-140\\
4 & 128s^4 - 128 s^3 + 1360s^2 - 976s + 612 & 256s^4 - 256s^3 + 4256 s^2 - 2720s + 4200\\
 \hline
\end{array}
$$
\caption{Values of $p_n(s)$ and $q_n(s)$.}
\label{ta71}
\end{table}

From \eqn{414} and \eqn{415} we deduce that
that these polynomials  satisfy
the  three-term recurrence relations
\beql{416}
p_{n+1} (s) = s p_n (s+2) + (2s-1) p_n (s) + (s-1) p_n (s-2)\,,
\eeq
and
\beql{417}
q_{n+1} (s) = (s+ 1) q_n (s+2) + (2s-1) q_n (s) + (s-2) q_n (s-2) \,.
\eeq
Three-term recurrence relations are associated to
orthogonal polynomials, with respect to a suitable measure, which
is explicitly 
given in the following result.

\begin{prop}\label{nle83} (Orthogonal polynomials)
For $n \ge 0$ the  polynomials
\beql{418aa}
P_n(x) := p_n (\frac{1}{2} + ix)
\eeq
are orthogonal polynomials with respect to the measure
$| \Gamma ( \frac{1}{4} + \frac{ix}{2} ) |^2 dx$ on the real $x$-axis.
For $n \ge 0$, the   polynomials
\beql{418bb}
Q_n(x) :=q_n (\frac{1}{2} + ix)
\eeq
are orthogonal polynomials with respect to the measure
$| \Gamma ( \frac{3}{4} + \frac{ix}{2} ) |^2 dx$ on the real $x$-axis.
\end{prop}

\pf
This was  observed by
Bump et~al.\cite[pp. 3-- 4]{BCKV99}. These polynomials are certain Meixner-Pollaczek
polynomials  $P_n^{(\lambda)}(x; \theta)$, in
the notation of  Koelink and Swarttouw \cite[Sect. 1.7]{KS98}. Namely, 
$$
P_n(x) = P_n^{(\frac{1}{4})}(\frac{x}{2}; \frac{\pi}{2}) :=
\frac{(\frac{1}{2})_n}{n!} i^n~
 {}_2F_1( -n, \frac{1}{4}+ \frac{ix}{2}; \frac{1}{2}; 2),
$$
and
$$
Q_n(x) = P_n^{(\frac{3}{4})}(\frac{x}{2}; \frac{\pi}{2}) :=
\frac{(\frac{3}{2})_n}{(n!}i^n 
~{}_2F_1( -n, \frac{3}{4}+ \frac{ix}{2}; \frac{3}{2}; 2).
$$
Here  the rising factorial $(\lambda)_n := \lambda(\lambda+1)\cdots (\lambda+n-1)$,
see also
Kutznetsov \cite{Ku07}, \cite{Ku08}.
 ~~~$\bsq$ \\

The Mellin transform and the multiplicative averaging operator $A^{a, c}$  given
in definition \ref{nde20}
map the infinite
dimensional vector space $\sS_0$ to a vector space $V_{a,c}$
of functions in the $s$-variable, via 
\beql{750}
\phi(x) \mapsto \sM( A^{a,c}[\phi])(s).
\eeq
This map is an intertwining map defining a Weyl algebra action
on the vector space $V_{a, c}$.  The image of the
Hermite basis under this intertwining gives
a basis of  the vector space $V_{a,c}$,
which we show below is expressible  in terms of Lerch zeta functions.


\begin{defi}~\label{de81}
{\em For each $n \in \ZZ_{\ge 0}$ and fixed 
$(a, c) \in \Box^{\circ}$ the {\em generalized (completed) Lerch zeta
function} $\hzt_n(s, a,c)$ is }
\begin{eqnarray}\label{418a}
 \hzt_n (s,a,c)& := & \frac{1}{(2\pi )^{n/2}} \sM(2A^{a,c}[\phi_n])(s) \nonumber \\
& = &\frac{2}{(2 \pi )^{n/2}} \int_{0}^{\infty} \left(\sum_{m \in \ZZ}
\phi_n( (m+c)x)e^{2\pi i m a} \right)x^{s-1}dx.
\end{eqnarray} 
\end{defi}

Note that for  $n=0, 1$ the function $\hzt_0(s,a,c)= \hat{L}^+(s, a,c)$
and $\hzt_1(s, a,c) = \hat{L}^{-}(s, a,c)$.
The  generalized Lerch zeta functions $\hzt_n (s,a, c)$
defined using the Hermite basis functions $\phi_n$ have a simple relation to
the  functions $L^{\pm}(s, a, c)$,  and satisfy similar
functional equations, as follows.
\begin{theorem}\label{th81}
For each $n \in \ZZ_{\ge 0}$ and fixed real $(a, c)$ with $0< a, c < 1$
the generalized (completed) Lerch zeta function
 \beql{418} 
 \hzt_n (s,a,c) :=  \frac{1}{(2\pi )^{n/2}}
\sM(2A^{a,c}[\phi_n])(s) = \frac{2}{(2 \pi)^{n/2}} \int_0^\infty A^{a,c}[\phi_n](x)~ x^{s-1} dx
 \eeq
extend to entire functions of $s$.  These functions are given by
\begin{eqnarray}\label{421}
\hzt_{2n} (s,a,c) & = &
\frac{1}{(2 \pi)^{n/2}} p_n(s) \hat{L}^{+}( s,a,c),  \\
\label{422}
\hzt_{2n+1} (s,a,c) & = &
\frac{1}{(2 \pi)^{n/2}} \sqrt{2\pi}\, q_n(s) \hat{L}^{-}( s,a,c) ~.
\end{eqnarray}
For each $n \ge 0$ they satisfy the   functional equation
\beql{419}
\hzt_n (s,a,c) = i^n~ e^{- 2 \pi iac} \hzt_n (1-s, 1-c, a) ~.
\eeq
\end{theorem}

\pf
It follows from the definition of $\hzt_n (s,a,c)$ and the
relation $2A^{a,c} = A_0^{a,c} + A_1^{a,c}$, together with Lemma \ref{nle21}, that 
\beql{420}
 \hzt_n (s,a,c) = \frac{1}{(2 \pi )^{n/2}} \{ \sM_0
(\phi_n)(s)L^{+}( s,a,c) +
 \sM_1 (\phi_n)(s)L^{-}(s,a,c) \} ~.
\eeq
Now $\phi_{2n} (x)$ is an even function and
$\phi_{2n+1} (x)$ is an odd function, so that
$\sM_k (\phi_n)(s) =0$ if $k \not\equiv n$ $(\bmod~2)$. Thus we obtain
\eqn{421} and \eqn{422} according as $n$ is even or odd.
The functional equations \eqn{419} now follow from that of
$\hat{L}^{\pm}( s, a,c)$  in Theorem \ref{th11} and of
$p_n (s)$ and $q_n (s)$ in Lemma \ref{nle42}.
The entire function property of each $\hzt_n(s,a,c)$ in the $s$-variable is
inherited from $\hat{L}^{\pm} (s,a,c)$   in Theorem~\ref{th11} and the fact that
$p_n (s)$ and $q_n (s)$ are polynomials.
 ~~~$\bsq$ \\

\paragraph{\bf Remarks.}
 (1) Theorem~\ref{th81} shows that the action of the Fourier transform on $\sS_0$ 
carries over 
an operator of order $4$ acting the vector space $V_{a,c}$,
which is compatible with the Weyl algebra action.

(2) The results of Bump et~al.\cite{BCKV99}
for the families of polynomials $\{p_n(s)\}$ (resp. $\{q_n(s)\}$)
may alternatively be
viewed as attached to the Lie algebra ${\bf sl}(2)$, viewed as
\beql{Z504}
{\bf sl}(2) = \RR \left[
x^2, \frac{1}{2}
\left( x \frac{\partial}{\partial x} + \frac{\partial}{\partial x} x \right),
\frac{\partial^2}{\partial x^2} \right],
\eeq
rather to its universal enveloping algebra $U( {\bf sl} (2))$.
The universal enveloping algebra $U( {\bf sl} (2))$ is a
subalgebra of the Weyl algebra  $\bA_1=U({\bf h}_{\RR})$.
As  a vector space over $\RR$, $U( {\bf sl} (2))$ is spanned by all
monomials in $x$ and $\frac{\partial}{\partial x}$
of even degree in $\bA_1$, counting $x$ and
$\frac{\partial}{\partial x}$ as degree one.
The members
of each family  $\{p_n(s)\}$ (resp. $\{q_n(s)\}$) (multiplied by
scalar constants)  are produced from the
bottom element $p_0(s)$ (resp. $q_0(s)$ by repeatedly applying
the raising operator $D_{+}^2$.

(3) For comparison with the real Heisenberg Lie
algebra treated here, the complex Heisenberg algebra
$\bh_{\CC} = \bh_{\RR} \otimes_{\RR} \CC$ used in quantum mechanics
has position operator $p=M_x$ and momentum operator 
$q = \frac{1}{i} \frac{\partial}{\partial x}$, so that $[p,q] = -iI$. The creation
operator $a = p + \frac {i}{2\pi} q$ and annihilation operator
 $a^*= p - \frac{i}{2 \pi} q$ satisfy 
 $[a, a^*] = I$ and $H= \frac{1}{2}(a a^* + a^*a).$

%
%
%
%
%


\end{document}

%% file: fg1.pstex_t
\begin{picture}(0,0)%
\special{psfile=fg1.pstex}%
\end{picture}%
\setlength{\unitlength}{2960sp}%
\begingroup\makeatletter\ifx\SetFigFont\undefined
\def\x#1#2#3#4#5#6#7\relax{\def\x{#1#2#3#4#5#6}}%
\expandafter\x\fmtname xxxxxx\relax \def\y{splain}%
\ifx\x\y   
\gdef\SetFigFont#1#2#3{%
  \ifnum #1<17\tiny\else \ifnum #1<20\small\else
  \ifnum #1<24\normalsize\else \ifnum #1<29\large\else
  \ifnum #1<34\Large\else \ifnum #1<41\LARGE\else
     \huge\fi\fi\fi\fi\fi\fi
  \csname #3\endcsname}%
\else
\gdef\SetFigFont#1#2#3{\begingroup
  \count@#1\relax \ifnum 25<\count@\count@25\fi
  \def\x{\endgroup\@setsize\SetFigFont{#2pt}}%
  \expandafter\x
    \csname \romannumeral\the\count@ pt\expandafter\endcsname
    \csname @\romannumeral\the\count@ pt\endcsname
  \csname #3\endcsname}%
\fi
\fi\endgroup
\begin{picture}(8871,2046)(774,-3494)
\put(8701,-2986){\makebox(0,0)[b]{\smash{\SetFigFont{9}{10.8}{rm}$c$}}}
\put(8701,-3436){\makebox(0,0)[b]{\smash{\SetFigFont{9}{10.8}{rm}$\Re (s) > 1$}}}
\put(5401,-2986){\makebox(0,0)[b]{\smash{\SetFigFont{9}{10.8}{rm}$c$}}}
\put(5401,-3436){\makebox(0,0)[b]{\smash{\SetFigFont{9}{10.8}{rm}$0 \le \Re (s) \le 1$}}}
\put(4351,-2236){\makebox(0,0)[rb]{\smash{\SetFigFont{9}{10.8}{rm}$a$}}}
\put(1051,-2236){\makebox(0,0)[rb]{\smash{\SetFigFont{9}{10.8}{rm}$a$}}}
\put(7651,-2236){\makebox(0,0)[rb]{\smash{\SetFigFont{9}{10.8}{rm}$a$}}}
\put(2101,-3436){\makebox(0,0)[b]{\smash{\SetFigFont{9}{10.8}{rm}$\Re (s) < 0$}}}
\put(2101,-3061){\makebox(0,0)[b]{\smash{\SetFigFont{9}{10.8}{rm}$c$}}}
\end{picture}